\documentclass[reqno,oneside, notitlepage]{amsart}

\usepackage{geometry}
\usepackage[dvipsnames]{xcolor}
\numberwithin{equation}{section}

\usepackage{tikz}
\usetikzlibrary{arrows.meta, positioning}
\usetikzlibrary{arrows.meta, automata,  
                positioning,
                quotes}

\usepackage{amsmath}
\usepackage{amssymb, amsthm}
\usepackage{mathtools}
\usepackage{mathrsfs}
\usepackage{esint}
\usepackage{bm}

\usepackage{graphicx}
\usepackage{tikz}
\usepackage{subcaption}
\usetikzlibrary{patterns} 
\usepackage{tikz-3dplot}
\tikzset{tdplot_fillconestyle/.style={fill=gray!10,opacity=0.3}}
\tdplotsetmaincoords{70}{25}

\usepackage[backend=biber,style=numeric,maxnames=99,minnames=1]{biblatex}
\AtEveryBibitem{%
  \clearfield{url}%
  \clearfield{issn}%
  \clearfield{isbn}%
} 
\DeclareFieldFormat[article]{volume}{\mkbibbold{#1}}  
\DeclareFieldFormat[article]{number}{\mkbibparens{#1}}  

\renewbibmacro*{volume+number+eid}{
  \printfield{volume}%
  \setunit{}%
  \printfield{number}%
  \setunit{\addcomma\space}%
  \printfield{eid}}

\AtBeginBibliography{\small}

\usepackage{stackengine}
\usepackage{IEEEtrantools}

\newtheorem{theorem}{Theorem}[section]
\newtheorem{lemma}[theorem]{Lemma}
\newtheorem{proposition}[theorem]{Proposition}

\theoremstyle{definition}

\theoremstyle{remark}
 \newtheorem{remark}[theorem]{Remark}

 \newcommand{\Energy}[1]{E_{\varepsilon}^{(#1)}(u;\chi)}
 \newcommand{\Energyinf}[1]{E_{\varepsilon}^{(#1)}(\mathcal{A}_F)}

 \newcommand\subsetsim{\mathrel{%
  \ooalign{\raise0.2ex\hbox{$\subset$}\cr\hidewidth\raise-0.8ex\hbox{\scalebox{0.9}{$\sim$}}\hidewidth\cr}}}

 \newcommand{\Eelinf}[1]{E_{\rm el}^{(#1)}  (\chi; F)}
 \newcommand{\Eel}[1]{E_{\rm el}^{(#1)} (u; \chi)}

\newcommand{\Esurf}{E_{\rm surf} (\chi)}
\DeclareMathOperator{\Per}{Per(\Omega)}
\newcommand{\diag}{{\rm diag}}

\DeclareMathOperator{\dist}{dist}

\newcommand{\dx}{\di x}

\newcommand{\N}{\mathbb{N}}
\newcommand{\Z}{\mathbb{Z}}
\newcommand{\R}{\mathbb{R}}

\newcommand{\defas}{\coloneqq} 
\newcommand*{\di}{\mathop{}\!\mathrm{d}}

\newcommand{\loc}{{\rm loc}}

\usepackage{todonotes}
\usepackage{marginnote}

\DeclarePairedDelimiterX\setof[1]\{\}{#1}
\DeclarePairedDelimiterX\abs[1]\lvert\rvert{#1}
\DeclarePairedDelimiterX\norm[1]\lVert\rVert{#1}
\DeclarePairedDelimiterX\sprod[2]\langle\rangle{#1, #2}

\usepackage{enumitem}

\setlist[enumerate,1]{font=\normalfont}
\setlist[itemize,1]{font=\normalfont}
 
\newlist{thmlist}{enumerate}{1}
\setlist[thmlist]{label=(\roman{thmlisti}),
	ref=(\roman{thmlisti}),font=\normalfont,
	noitemsep}

\usepackage[colorlinks=true, allcolors=blue]{hyperref} 
\usepackage[toc]{appendix}
\usepackage{url}
\usepackage{xurl} 

\usepackage{color}

	\title{Sharper energy scaling bounds of the singularly perturbed Tartar square}%

		\title[Sharper energy scaling bounds of the singularly perturbed Tartar square]{Sharper energy scaling bounds of the singularly perturbed Tartar square}

		\author{Lennart Machill} 
\address{Institute for Applied Mathematics, University of Bonn, Endenicher Allee 60, 53115 Bonn, Germany}
\email{lmachill@uni-bonn.de}
	
\begin{document}

	\begin{abstract}
		In this article, we discuss the scaling behaviour of the singularly perturbed Tartar square, revisiting the results by {\sc Rüland and Tribuzio} \cite{RT22}. 
Therein, it was proven that both lower and upper bounds satisfy a behaviour of the same type, up to a small mismatch of logarithmic powers.
We show by an improvement of the ansatz free lower bound that the logarithmic powers indeed match and that the bounds are of order $\exp( - C \log(\varepsilon^{-1})^{1/2} )$ 
for suitable constants $C = C_{\rm low}$ and $C = C_{\rm upp}$, where $\varepsilon$ denotes the singular perturbation parameter.
 More specifically, the bootstrap argument in the Fourier space is refined by resorting to Fourier multipliers enjoying $L^\infty$-bounds which are constructed via de la Vallée Poussin kernels.
Regarding upper bounds, we perform iterated branching constructions sharpening the upper bound of {\sc Winter} \cite{W97} which resorted to an iteration of laminates as a competitor.
Eventually, the mismatch of upper and lower bound is quantified, and for quadratic energies we obtain the bounds $C_{\rm low} \leq 7.0339$ and  $C_{\rm upp} \geq 2 \log(2) - \nu $ for any $\nu>0$. 
\end{abstract}


	\maketitle

\setcounter{tocdepth}{2}	

\section{Introduction}
Given a bounded open set $\Omega \subset \R^d$, differential inclusions of the form
\begin{align}\label{eq:differentialinclusion}
	\nabla u(x) \in K \quad \text{for almost every } x \in \Omega
\end{align}
and functions $u \colon \Omega \to \R^d$
have important applications in the calculus of variations and arise e.g.\ in the mathematical modeling of elastic materials.
Here, $K \subset \R^{d \times d}$ may represent different states or stress-free configurations of a solid.
A particularly interesting class of such materials are shape-memory alloys which undergo a solid-to-solid, diffusionless phase transformation in the passage from the high to the low temperature phase. In the latter, crystalline symmetry is reduced, leading to a highly non-convex set $K$ with multiple preferred states \cite{B04,B}.

Here, both the structure of $K$ and the regularity of $u$ strongly affects the class of solutions.
With the aim of modelling and studying the emergence of complex patterns and microstructures in shape-memory alloys, it is interesting to understand in which cases the inclusion \eqref{eq:differentialinclusion} is \emph{flexible}, corresponding to the presence of non-trivial solutions (cf.\ \cite{M1}).
Classically, non-trivial piecewise affine and continuous solutions exist if there are two state $A \neq B \in K$ with  ${\rm rank}(A-B) =1$, the so-called twins or more generally \emph{laminates} (see e.g.\ \cite{BallJames87}).
Surprisingly, the lack of rank-one connections in $K$ does not characterize \emph{rigid} differential inclusions, namely those for which \eqref{eq:differentialinclusion} has only affine solutions.
Indeed, Kirchheim and Preiss \cite{Kirchheim2003} 
provided a striking example of a flexible discrete set $K$ of \emph{five} elements via convex integration (see also \cite{ForsterSzekelyhidi18}).
Remarkably, five is indeed the critical cardinality for a set $K$ without rank-one connections to be flexible, while Lipschitz solutions to \eqref{eq:differentialinclusion} are rigid 
if $K$ consists less than five elements \cite{chlebikkirchheim2002,Sverak1991}.

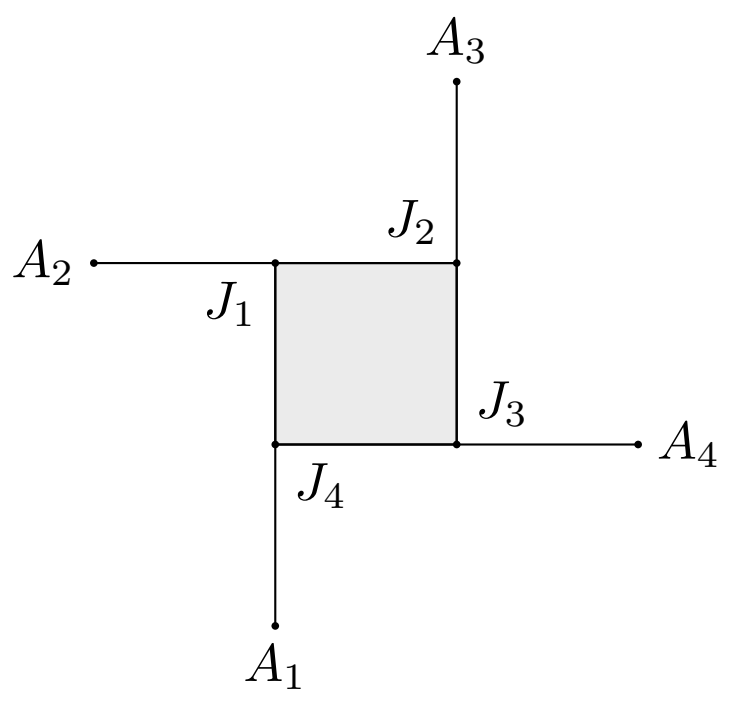
\begin{figure}[t]
\begin{tikzpicture}[scale=0.6]
\coordinate (J1) at (-1, 1);
\coordinate (J2) at ( 1, 1);
\coordinate (J3) at ( 1,-1);
\coordinate (J4) at (-1,-1);

\coordinate (A1) at (-1,-3);
\coordinate (A2) at (-3, 1);
\coordinate (A3) at ( 1, 3);
\coordinate (A4) at ( 3,-1);

\filldraw[fill=black, fill opacity=0.08, draw=black]
    (J1) -- (J2) -- (J3) -- (J4) -- cycle;

\draw (A2) -- (J2);
\draw (A3) -- (J3);
\draw (A4) -- (J4);
\draw (A1) -- (J1);

\fill (J1) circle (1.2pt) node[below left]  {$J_1$};
\fill (J2) circle (1.2pt) node[above left] {$J_2$};
\fill (J3) circle (1.2pt) node[above right] {$J_3$};
\fill (J4) circle (1.2pt) node[below right]  {$J_4$};

\fill (A1) circle (1.2pt) node[below] {$A_1$};
\fill (A2) circle (1.2pt) node[left]  {$A_2$};
\fill (A3) circle (1.2pt) node[above] {$A_3$};
\fill (A4) circle (1.2pt) node[right] {$A_4$};

\end{tikzpicture}
	\caption{Visualization of the Tartar square. The union of the lines and the filled square represent the rank-one convex, quasiconvex and polyconvex hulls $K^{\rm rc} = K^{\rm qc} = K^{\rm pc}$ respectively. The matrices $J_i$, $i = 1, \dots,4 $, are defined in \eqref{eq:auxmatrices}, and the rank-one directions correspond to horizontal and vertical lines.}
	\label{fig:Tartarsquare}
\end{figure}

Such a dichotomy between rigidity and flexibility is also visible on the level of approximate solutions, i.e., 
 for sequences $(u_j)_{j \in \N}$ with $u_j \in W^{1,\infty}(\Omega,\R^2)$ such that $\dist(\nabla u_j, K) \to 0$ in measure.
In this case, the critical cardinality has proved to be \emph{four}.
If $K$ has less than four elements which are not rank-one connected,
every $(\nabla u_j)_{j \in \N}$ converges to an affine solution \cite{FaracoSzekelyhidi2008}. 
Contrary, the problem can become flexible if $K$ consists of four matrices.
The so-called Tartar square \cite{Tartar93}, whose wells are given by
	\begin{align}\label{eq:choiceofwells}
		K  \defas \{ A_1, \dots, A_4\} \defas \{ \diag (-1,-3),  \ \diag(-3,1) , \ \diag(1,3) , \ \diag(3,-1) \},
	\end{align}
	displays such a flexibility.
	 Although it is rigid on the level of exact solutions \cite{chlebikkirchheim2002},
 it is a $T_4$-structure, i.e., the wells do not share rank-one connections, but admit cyclic rank-one splittings. This allows to construct approximate solutions to the associated differential inclusions 
 which, however, \emph{do not} converge to an affine solution.
The associated semiconvex hulls are visualized in Figure~\ref{fig:Tartarsquare}.  
We refer to \cite{Ri18,M1} for a more detailed introduction on rigidity and flexibility.
	
As in \cite{RT22,W97}, we seek to investigate the microstructure of approximate solutions for the Tartar square quantitatively by analysing a singular perturbation of its associated elastic energy. 
To this end, let us fix $\Omega \defas (0,1)^2 \subset \mathbb{R}^2$ and analyse the scaling behavior of the singularly perturbed problem
\begin{align}
\label{eq:minprob}
	\Energyinf{p}  \defas \inf\limits_{\chi \in BV((0,1)^2; K)} \inf\limits_{u \in \mathcal{A}_F  } \Energy{p},
\end{align}
for sufficiently small values of $\varepsilon>0$, where we set
\begin{align}\label{def:periodicboundarycond}
\mathcal{A}_F \defas \{    u \in W^{1,p}_{\loc}(\R^2; \R^2) : \langle \nabla u \rangle = F , \  \nabla    u \text{ is } [0,1]^2\text{-periodic} \}            
\end{align}
and
	\begin{align}
		\Energy{p} & \defas \Eel{p} + \varepsilon \Esurf  \defas \int_{(0,1)^2}\vert \nabla u - \chi \vert^p \dx + \varepsilon  \Vert D \chi \Vert_{TV((0,1)^2)}. \label{def:phaseindic}
	\end{align}  
After upper scaling bounds for the singularly perturbed Tartar square were investigated in \cite{W97} and \cite[Section~4]{Chipot99}, an ansatz free lower bound was derived only recently in \cite{RT22}.
Our main goal is to sharpen the latter result by showing the following improved lower energy estimate.
\begin{theorem}[Lower bound]\label{mainth}
Let $p \in [2,+\infty)$ and $F \in K^{\rm qc}\setminus K$.
 Then, there exists a constant $c_{\rm low} = c_{\rm low}(K,p, F)>0$ and a constant $C_{\rm low} = C_{\rm low}(K,p)>0  $ 
 such that for every $\varepsilon \in (0,1)$ it holds that
\begin{align*}
	c_{\rm low} \exp\left(- C_{\rm low}  \log(\varepsilon^{-1})^{\frac{1}{2}} \right) \leq \Energyinf{p}.
\end{align*}  
Further,  the best constant $C_{\rm low}$ satisfies at least the upper bound $C_{\rm low} \leq  \max\{ 4.9730 \sqrt{p}, 3.2243 p\}$.
\end{theorem}
Let us briefly compare the bound with the result of \cite{RT22}.
Therein, a lower scaling bound of the form
$ c_{\rm low,\nu} \exp(-C_{\rm low,\nu} \log(\varepsilon^{-1})^{\frac{1}{2} + \nu}) \leq \Energyinf{2}$ has been shown
for every $\nu>0$. The constant $C_{\rm low,\nu} $, however, degenerates as $\nu \to 0$. Our result hence sharpens the bound and shows its validity for $\nu = 0$.
Further, we also provide lower scaling bounds for $p\geq 2$ in a similar manner to \cite{LMAR26} by exploiting that $L^p$-based elastic energies control the $L^2$-based counterpart in this case. 
The explicit upper bound on $C_{\rm low}$ result is strongly influenced by technical obstructions in the Fourier-based framework. Its sharpness is discussed after the proof in Remark~\ref{rem:optimalconst}.
Here, we note that the scaling behaviour persists for the problem 
  $$\inf\limits_{u \in A_F, \nabla u \in BV} \int_{(0,1)^2}\dist^p(\nabla u , K) \di x + \varepsilon \Vert D^2 u \Vert_{TV((0,1)^2)},$$
  which is proved in \cite{ARCTATCZ25}. The advantage of the formulation via phase indicators renders the problem easier for a Fourier-based approach as the elastic energy becomes convex in $u$.
Regarding upper bounds, we show the following:
\begin{theorem}[Upper bound]\label{mainth:upper}
Let $p \in [1,+\infty)$ and $F \in K^{\rm qc}\setminus K$.
 Then, there exists a constant $c_{\rm upp} = c_{\rm upp}(K,p,F)>0$ and a constant $C_{\rm upp} = C_{\rm upp}(K,p)>0  $ such that for every $\varepsilon \in (0,1)$ it holds that
\begin{align*}
	 \Energyinf{p} \leq c_{\rm upp} \exp\left(-C_{\rm upp}  \log(\varepsilon^{-1})^{\frac{1}{2}} \right) .
\end{align*}  
Further,  the best constant $C_{\rm upp}$ satisfies at least the lower bound $C_{\rm upp} \geq \sqrt{2p \log(2)} - \nu $ for any $\nu>0$   . 
\end{theorem} 
Note that our upper bound construction does not allow for the choice $C_{\rm upp} =  \sqrt{2p \log(2)}$ as $c_{\rm upp} \to \infty$ as $C_{\rm upp} \to \sqrt{2p \log(2)}$.
Such a degeneracy has already been observed in \cite{W97} for iterated laminate constructions where it was shown that $C_{\rm upp} > \sqrt{2 \log(2)} -\nu$ for quadratic elastic energies. 
In contrast, our upper bounds resorts to iterated branching constructions, showing a faster scaling bound for $p>1$.
Note that for $L^1$-based elastic energies, branching and laminate constructions lead to the same scaling behaviour which is well-known for boundary data of finite order (see \cite{Lorent2001,ARCTATCZ25} and the discussion in \cite{RT23}).
In contrast to data of finite order \cite{LMAR26}, the in Theorem~\ref{mainth} obtained lower bound for a general $p$ does not coincide with the behaviour in $p$ in Theorem~\ref{mainth:upper}.
Further, let us stress that the constants in Theorem~\ref{mainth} and Theorem~\ref{mainth:upper} do not match, but our method significantly improves the constants from \cite{RT22} and particularly allows for an explicit upper bound of $C_{\rm low}$. Hence, we are able to quantify the gap between upper and lower bound.

As discussed in \cite{W97,RT22}, the scaling behaviour is slower than a polynomial scaling, but faster than a logarithmic scaling, in the sense that
\begin{align*}
	\varepsilon^\gamma \ll \exp (-C \log(\varepsilon^{-1})^{1/2}) \ll \log(\varepsilon^{-1})^{-\gamma}
\end{align*}
for all $\gamma>0$, showing a high complexity of microstructure.
This particular scaling behaviour persists if one replaces the periodic boundary conditions in
 \eqref{def:periodicboundarycond}
with Dirichlet boundary conditions \cite{RT22}.
In fact, they are more restrictive and can be viewed as a smaller subset of the periodic setting in the minimization problem \eqref{eq:minprob} via the divergence theorem. 
This shows that our lower bound stays admissible while 
our upper bound construction is carried out with Dirichlet boundary conditions, i.e., the constants $C_{\rm low}$ and $C_{\rm upp}$ in Theorems \ref{mainth} and \ref{mainth:upper} remain valid for this modification.

Let us connect our article to further results in the literature.
Although Tartar initially conjectured that the differential inclusion \eqref{eq:differentialinclusion} is rigid for approximate solutions whenever the non-empty compact set $K$ does not contain rank-one connections \cite{Tartar82},
he himself discovered it to be false in \cite{Tartar93} resorting to the wells in \eqref{eq:choiceofwells}.
For connected sets, the conjecture turns out to be true \cite{Sverak1993},
and by excluding $T_4$-structures, the differential inclusion regains the rigidity as well \cite{FaracoSzekelyhidi2008}.
At this point, we also mention that the Tartar square was discovered in different contexts, see e.g.\ \cite{AumannHart1986,Casadio93,NesiMilton91,Scheffer74,Bhatetal94}.



 Regarding energy scaling laws, let us particularly highlight the seminal works by {\sc Kohn and Müller} \cite{KM1, KM2} which for the first time rigorously proved that 
 branching structures can be energetically preferable to simple laminate structures in terms of energy scaling, recovering $\varepsilon^{2/3}$-bounds for quadratic elastic energies.
In this context, we refer to
\cite{CC14,CC15} for scaling laws in
 the geometrically nonlinear framework, to \cite{RRT23,RRTT24,LMAR26} for a further systematic study of the influence of gauge invariances on the scaling law, and to \cite{contimisiatsupper2026} for upper bound constructions in more general geometries.
 The articles \cite{PF25,MRP26} deal with scaling laws associated with incompatible wells, and
 we further highlight the works \cite{GZ23,GRTZ25} on scaling laws for low energy microstructures, 
the investigation of nucleation results \cite{KK11,CDMZ20,AKKR24,RT23a,TZ24,KKO13,KO19} as well as models inspired by the study of hysteresis \cite{Z14,CDZ17, CZ16}. 
 Scaling laws are also investigated in the context of compliance minimization \cite{KW14, KW16, PW22}, pattern formation in biomembranes \cite{Anastasija26}, and micromagnetism \cite{ChoksiKohnOtto99,KnupferMuratovNolte19,KnupferStantejsky22,GinsterZwicknagl24,RiedRoman25}.
Further, we mention  \cite{BellaMarziani25,Marziani2026} for recent developments on wrinkling phenomena.
  Different forms of surface energy regularizations are content of \cite{ARCTATCZ25} and finer properties of minimizers for energies of order $\varepsilon^{2/3}$ have been studied in \cite{C1, CDKZ21, TS21, TS21a}. 
We also refer to a recent scaling law \cite{CoKoOl24}, analysing the relation of bending effects with twinning.

Triggered by the rigidity estimates for the
geometrically linearised cubic-to-tetragonal phase transformation model  \cite{CO1}, Fourier analytic approaches have proved to be a successful method when dealing with scaling laws, see e.g.\ the previously mentioned compliance minimization \cite{KW14, KW16} or nucleation problems \cite{KKO13,KO19}.
They are robust for models including higher-order differential operators \cite{RRTT24} and for two-state problems including Dirichlet boundary conditions \cite{RRT23,PF25}.
Inspired by a model in compliance minimization \cite{KW16}, they also allow for the treatment of more complicated boundary data:
firstly, for sets of wells involving higher-order laminates of finite order \cite{RT23,RT23a,LMAR26}, and, secondly, for certain $T_N$-structures.
Building on the quantitative analysis of the Tartar square in \cite{RT22}, this $T_4$-structure was also analysed in a discrete model \cite{ARCTATCZ25}.
Although $T_3$-structures do not exist for differential inclusions on the gradient, they do for
the divergence and the symmetric gradient as they are more flexible operators. We refer to \cite{RRT23,IKRTZ26} for corresponding quantitative results.

In all the previously mentioned results on the quantitative analysis of $T_N$-structures the bootstrap argument in the Fourier space generated a logarithmic \emph{loss} in the lower bound.
The purpose of this article is to get rid of this obstacle. Our main observation is that the method of \cite{RT22} can also be employed on stripes instead of finite cones which allows for a larger pool of available Fourier multipliers in order to truncate high frequencies.
In particular, by resorting to de la Vallée Poussin kernels that enjoy $L^\infty$-estimates, technical estimates can be improved and quantified, see Proposition~\ref{prop:rectanglereduction} in Subsection~\ref{sec:Fourier5} below.
Hence, our result not only improves the scaling behaviour for $T_N$-structures,
it also sharpens the constant of lower bounds for general higher-order data that are derived via localisation in Fourier space.

Let us also mention that Fourier localisation methods are only applied for elastic energies that are at least quadratic,
while real-space localisation methods have already been successfully employed for a general $p \in [1,+\infty)$, see e.g.\ \cite{CDZ17} or \cite[Section~3.2]{RT23}.
The latter results address first-order boundary data while results for more complicated boundary data without tools from Fourier analysis are scarce.
We refer to \cite{Lorent2001} for a second-order boundary datum for $p =1$ in a discrete model, which was lifted only recently to a general $p \in [1,+\infty)$ in \cite{HoMaRue26}.

Finally, the singular perturbation problem \eqref{def:phaseindic} is closely related to the study of classical Modica--Mortola energies \cite{ModicaMortola77},
and their derivation of sharp-interface limits. In this context, one considers diffusive perturbations and energy scalings of order $\varepsilon$ which is in stark contrast to the bounds in Theorem~\ref{mainth} and Theorem~\ref{mainth:upper}.
We refer to \cite{ContiFonsecaLeoni02,ContiSchweizer06lin,ContiSchweizer06nonlin,KitavtsevLuckhausRueland17,FriedrichDavoli25} 
for solid-to-solid transformation settings, to \cite{SavinValdinoci12,Heilmann25,HeilmannFusco26,Gianna26,Davide25} for more general perturbations, and to \cite{Steinke2025,Kreutz1,Kreutz2} for second-order $\Gamma$-limits.

\section{Lower bound}
\label{sec:lowerbound}
   Our goal is to prove the lower bound
via a Fourier-based framework. 
 For a $\mathbb{T}^d$-periodic function $f \in L^2\left(\mathbb{T}^d\right)$ the Fourier transform at $\xi \in \mathbb{Z}^d$ is given by 
 $$
\mathcal{F}(f)(\xi) = \int_{\mathbb{T}^d} e^{- 2 \pi i \xi \cdot x} f(x) \, {\rm d}x ,
$$
 where the evaluations at $\xi$ are the Fourier coefficients.
The general strategy is to control the $L^2$-norm of the phase indicator $\chi$ in its representation as a Fourier series via elastic and surface energy.
More specifically, as $\chi$ does not attain the value $F$, we find a constant $c>0$ such that Plancherel's formula yields
	\begin{align}\label{eq:strategyFourier}
		 c \leq \Vert  \chi - F \Vert_{L^2}^2 = \sum\limits_{ \xi \in \mathbb{Z}^2 }  \vert \mathcal{F}(\chi - F)(\xi) \vert^2,
	\end{align}
	where one views $ 
		\chi - F $  as a $\mathbb{T}^d$-periodic function.
Then, one seeks to control the
Fourier coefficients on the right-hand side via elastic and surface energy in a quantitative form.
As these bounds will hold uniformly among the class of phase indicators,
one obtains a lower bound.

In Section~\ref{sec:Fourier1} and Section~\ref{sec:Fouriersurface}, we recall how the Fourier coefficients 
are controlled with elastic energy and surface energy in a quantitative way.
As these estimates do not allow for optimal scaling bounds,
the refined bootstrap argument is presented in Section~\ref{sec:Fourier5} and applied in Section~\ref{sec:proofslower}.

\subsection{The elastic energy in Fourier space} \label{sec:Fourier1}
 Let us reformulate the elastic energy as a Fourier series. 
 The following lemma controls Fourier coefficients in terms of elastic energy.
To this end, it is convenient to introduce the notation
$\hat \xi = \vert \xi \vert^{-1} \xi $ for normalized frequencies and the function
	$\tilde \chi =  
		\chi - F  $ viewed as a $ \mathbb{T}^d$-periodic function.
Moreover, we set
\begin{align*} 
	 \Eelinf{p} \defas \inf\limits_{u \in \mathcal{A}_F } \Eel{p} .
\end{align*}
\begin{lemma}[Control via elastic energy]\label{lem:Fourierchar}  
 Let $p \in [2,+\infty)$.  Then, it holds that
  	\begin{align*}  
	\Eelinf{p}^{2/p}   & \geq \vert \mathcal{F}( \tilde \chi) (0)  \vert^2 +  \sum\limits_{ \xi \in \mathbb{Z}^2 \setminus \{0 \}}   
		\vert \hat \xi_2 \vert^2 \vert \mathcal{F}(\tilde \chi_{11})(\xi) \vert^2 + \vert\hat \xi_1\vert^2 \vert \mathcal{F}( \tilde \chi_{22})(\xi) \vert^2 .
  \end{align*} 
  Further, for $i = 1,2$ and $\kappa>0$ we have
		\begin{align*}
		\sum_{\substack{\xi \in \Z^2\\ \dist(\hat \xi; \R e_i) > \kappa }}\vert \mathcal{F}(\tilde \chi_{ii})(\xi) \vert^2 \leq \kappa^{- 2 }   \Eelinf{p}^{2/p} .
	\end{align*}
\end{lemma} 
\begin{proof}
	The first inequality follows from \cite[Lemma~4.1]{RT23} in the case $p = 2$.
	The case $p>2$ is deduced from Hölder's inequality, see \cite[Appendix A.2]{LMAR26} for a rigorous treatment.
	The second bound is deduced from the first one as
$\dist^2(\hat \xi, \R e_1) = \vert \hat \xi_2 \vert^2 $ and $\dist^2(\hat \xi, \R e_2) = \vert \hat \xi_1 \vert^2 $.
\end{proof}
The result quantitatively shows for states with small elastic energy that the components $\mathcal{F}(\tilde \chi_{11})$ and $\mathcal{F}(\tilde \chi_{22})$ concentrate in the directions induced by the compatibility of the wells, see Figure~\ref{fig:Tartarsquare}.	
As the bound depends on \emph{normalized} frequencies, coercivity is obtained outside of \emph{cones} around the compatibility directions, which are precisely the directions $e_1$ and $e_2$.
 More specifically, jumps of the component $\tilde \chi_{ii}$ along a cone around the $e_i$-direction are energetically cheaper. 


   
 \subsection{The surface energy in Fourier space}\label{sec:Fouriersurface}
Fourier coefficients at high frequencies are controlled by the surface energy, as shown by the following lemma. Recall \eqref{def:phaseindic} for the definition of surface energy.
\begin{lemma}\label{lem:highfrequencyoutsidecones}
	 There exists a constant $C>0$ such that for all $\mu >0$ it holds that 
	\begin{align*}
		 \sum_{ \{ \xi \in \Z^2 : \, \vert \xi \vert \geq \mu \} } \vert \mathcal{F}( \tilde \chi_{ii} ) (\xi) \vert^2 \leq C \Vert \tilde \chi \Vert_{L^\infty} \mu^{-1} \left( \varepsilon^{-1} \Esurf +  \Vert \tilde \chi \Vert_{L^\infty} \Per \right).
	\end{align*} 
\end{lemma}
\begin{proof}
	The lemma is e.g.\ proved in \cite[Lemma~4.5]{RT23}, which itself is based on \cite[Lemma~4.3]{KKO13}.
	We also point to \cite[Lemma 2.2]{ARCTATCZ25} for a version with anisotropic surface energies, where the dependence on $\tilde \chi$ is made precise.
\end{proof}

\subsection{Quantitative interaction of Fourier coefficients}\label{sec:Fourier5}
As explained heuristically above, we aim at controlling the right-hand side of
\eqref{eq:strategyFourier}
in terms of elastic energy and surface energy, while the latter (in contrast to the former) generates the factor $\varepsilon^{-1}$ on the right-hand side, 
see Lemma~\ref{lem:Fourierchar} and Lemma~\ref{lem:highfrequencyoutsidecones}.
Hence, to obtain sharp lower bounds, as many Fourier coefficients as possible should be estimated purely with elastic energy.
In view of \eqref{eq:choiceofwells}, a computation shows 
$g(\chi_{11}) = \chi_{22}$ for
$g(t) = - \frac{5}{12} t^3 + \frac{41}{12} t$ as well as $-g(\chi_{22}) = \chi_{11}$. 
This means that the value $\chi_{22}$ is already determined by the value of $\chi_{11}$ and vice versa. 
This interaction will be made quantitative in terms of Fourier coefficients.
Roughly speaking, Fourier coefficients $\mathcal{F}(\tilde \chi_{11})(\xi_1,\xi_2)$ for large $\vert \xi_2 \vert$ are controlled by $\mathcal{F}(\tilde \chi_{22})$, 
which in turn may estimated by elastic energy as the cones point into a perpendicular direction. Hence, less surface energy is required.
The precise statement is addressed in Proposition~\ref{prop:rectanglereduction} below.

Let us stress once again that this general concept is not new and has already been employed for $T_N$-structures \cite{RT22,RRT23,ARCTATCZ25,IKRTZ26} as well as for sets of wells involving higher-order laminates \cite{RT23,LMAR26}.
We also highlight that this strategy has been employed successfully in the context of a nucleation problem \cite{RT23a}. 

This section shows that the estimates from the above-mentioned articles can be sharpened. For $T_N$-structures, this is particularly visible in scaling behaviour in $\varepsilon$. 
The main reason for the ``loss'' in these estimates are the lack of the bounds of the Fourier multipliers involved in these arguments:
although cones are naturally generated due to the control of the Fourier coefficients via elastic energy, see Lemma~\ref{lem:Fourierchar},
the Fourier multipliers which concentrate on cones cannot be controlled in the optimal way, eventually leading to the suboptimal lower bound discussed below Theorem~\ref{mainth}.
More precisely, given a cone
\begin{align*}
C_{\rm cone}^\gamma \defas \{\xi \in \Z^2: \dist(\hat \xi, V_j)  \leq C^{-1} \gamma, \ \vert \xi \vert \leq C \gamma \} \quad \text{ for } \quad \gamma>0,
\end{align*} 
the
associated Fourier multiplier $m_{\rm cone, \gamma}(D)$ satisfying $m_{\rm cone, \gamma}(\xi)=1$ for $\xi \in C_{\rm cone}^\gamma$ and $m_{\rm cone, \gamma}(\xi)=0$ for $\xi \in\Z^d \setminus C_{\rm cone}^{2\gamma}$ (as e.g.\ given in \cite[Section~4.3]{RT23})
satisfies the bound
\begin{align*}
	 \Vert m_{\rm cone, \gamma}(D) u \Vert_{L^q} \leq C \Vert u \Vert_{L^q} \text{ for } q \in (1,+\infty),
\end{align*}
but it is \emph{unclear} if the estimates transfer to $q = \infty$.
The bound for $q \in (1,+\infty)$ is based on a combination of the Marcinkiewicz’s multiplier theorem \cite[Corollary~6.2.5]{Grafakos} with the transference principle \cite[Theorem~4.3.7]{Grafakos}, see e.g.\ \cite[Proposition~2.1]{RT23}.
In essence, the lack of the control for $q = \infty$ is the main reason for the lower bound of the form $c_\nu\exp (- C_\nu \log(\varepsilon^{-1})^{1/2+ \nu})$ for $\nu>0$ in \cite{RT22}.

We overcome this technical difficulty by
replacing these multipliers with ones that are concentrated on \emph{stripes} and enjoy $L^\infty$-estimates with known constants.
More precisely, 
we consider multipliers generated by de la Vallée Poussin kernels. The latter correspond to a linear combination of Fejér kernels, which are uniformly bounded in $L^1$. Hence, the associated multiplier allows for an $L^\infty$-estimate, following from the convolution inequality, see \eqref{eq:conv} below.

\subsubsection*{De la Vallée Poussin kernels}
In view of \cite[Exercise~3.1.4]{Grafakos}, we consider de la Vallée Poussin kernels of the form
\begin{align}\label{eq:representatioPoussinkernel}
	V_{N}(x_1) = \begin{cases} 1& \quad \text{if } N = 0 , \\
	2 F_{2N +1}(x_1) - F_{N }(x_1)& \quad \text{if }N \in \N ,
	\end{cases}
\end{align} 
where $(F_N)_{N \in \N }$ are Fejér kernels. These kernels are averages of Dirichlet kernels which satisfy better bounds, namely $\Vert F_N \Vert_{L^1} \leq 1$, see \cite[Proposition~3.1.10]{Grafakos}.
Then, given a bounded periodic function $f \in L^\infty(\mathbb{T})$ and $N \geq 1$, a transformation of $V_N \ast f$ into Fourier space yields
\begin{align}\label{cut-offmultiplier}
	\mathcal{F}(V_N \ast f)(\xi_1) = \mathcal{F}(V_N)(\xi_1)\mathcal{F}(f) (\xi_1), \text{ where } \mathcal{F}(V_N)(\xi_1) = \begin{cases}
		1, & \vert \xi_1 \vert \leq N ,\\
		0, & \vert \xi_1 \vert \geq 2 N , \\
		T(\xi_1), & \vert \xi_1 \vert \in (N, 2N), \\
	\end{cases}
\end{align}
  where $T \colon \R \to [0,1]$ interpolates linearly between $1$ and $0$, while for $N=0$ only the frequency $\xi = 0$ is selected. Indeed, for the representation of $\mathcal{F}(V_N)$ we refer e.g.\ to \cite{MEHTA15} or \cite[Exercise~3.1.4]{Grafakos}.
Then,
  \eqref{eq:representatioPoussinkernel} and the $L^1$-bound on Fejér kernels imply by the convolution inequality
  \begin{align*}
	\Vert V_N \ast f \Vert_{L^\infty} \leq \Vert V_N \Vert_{L^1} \Vert f \Vert_{L^\infty} \leq 3 \Vert f \Vert_{L^\infty}
  \end{align*}
  for all $N \in \{0\} \cup \N$.
  It is worth noting that the presence of a Fejér kernel which is constructed by averages of Dirichlet kernels is crucial. 
An approach with 
Dirichlet kernels, which corresponds to a sharp truncation of frequencies, \emph{fails} as its $L^1$-norm blows up logarithmically, see e.g.\ \cite[Exercise~3.1.5]{Grafakos}.

  \subsubsection*{Improved bounds on the $L^1$-norm}
	As the $L^1$-norm of $V_N$ influences the constant $C_{\rm low}$ in our analysis below, one may wonder if the bound $3$ is optimal.
	This was determined in \cite[Corollary 1.3]{MEHTA15} showing that
	\begin{align} \label{eq:bestbound}
	\lambda \defas \Vert V_N \Vert_{L^1} = \frac{1}{3} + \frac{2 \sqrt{3}}{\pi} \approx 1.44.
	\end{align}
	 This constant is known to be the first Lebesgue constant which is defined by the $L^1$-norm of the first Dirichlet kernel.

  \subsubsection*{Definition of multipliers}
We now address the desired construction for our application.
  To this end, we define the stripes
\begin{align*}
	R_{i,\mu} \defas \{\xi \in \Z^2: \vert \xi_i \vert \leq \mu \text{ for } \vert \xi_i \vert^2 + \vert \xi_j \vert^2 = \vert \xi \vert^2 , \ j \neq i  \} \text{ for } i = 1,2, \ \mu > 0. 
\end{align*}
 We set
\begin{align*}
	\psi_{i,\mu}(D) f (x_i,x_j) \defas    V_{\lfloor\mu\rfloor}(x_i)   \ast f (x_i,x_j) 
\end{align*}
where $\lfloor\cdot \rfloor$ denotes the floor-function.
Recalling \eqref{cut-offmultiplier}, the multiplier $\psi_{i,\mu}(D)$, $i = 1,2$ satisfies $\mathcal{F}(\psi_{i,\mu})(\xi) = 1$ for $\vert\xi_i\vert \leq \lfloor\mu\rfloor$ and $\mathcal{F}(\psi_{i,\mu})(\xi) = 0$ for $ \vert \xi_i \vert \geq   2\lfloor\mu\rfloor  $.
Hence, the multiplier selects Fourier coefficients in a stripe around \emph{the $e_j$-axis}.
Moreover, the convolution inequality implies that
\begin{align}\label{eq:conv}
	\Vert \psi_{i,\mu}(D) \tilde \chi_{ii} \Vert_{L^\infty} \leq  \int_\mathbb{T}  \vert V_{\lfloor\mu\rfloor} (x_i) \vert   \di x_i   \Vert \tilde \chi_{ii} \Vert_{L^\infty}  \leq  \lambda \Vert \tilde \chi_{ii} \Vert_{L^\infty} ,
\end{align}
where $\lambda$ was defined in \eqref{eq:bestbound}.
The Fourier multipliers which are associated to Dirichlet kernels are denoted by $m_{i,\mu}(D)$ and satisfy $\mathcal{F}(m_{i,\mu})(\xi) = 1$ for $\vert\xi_i\vert \leq \lfloor\mu\rfloor$ and $\mathcal{F}(m_{i,\mu})(\xi) = 0$ for $ \vert \xi_i \vert >   \lfloor\mu\rfloor  $.
This corresponds to a sharp truncation.
With this in hand, we prove the following result.
\begin{proposition}[Stripe reduction]\label{prop:rectanglereduction}
Let $g$ be a polynomial with $\deg(g)\geq 1$ such that
 for $i \neq j$, $i, j =1,2$, it holds that
 $\tilde{\chi}_{jj} = g(\tilde\chi_{ii})$.
Consider $\kappa_i \in (0,1)$ and $\mu_i>0$ such that $\mu_j = 2 \deg(g) \mu_i \kappa_i (1-\kappa_i^2)^{-1/2}$.
 Then, it holds that   
\begin{align}\label{eq:estimatestripered}
 \|\tilde{\chi}_{jj} - m_{j,\mu_{j}}(D) \tilde{\chi}_{jj}\|_{L^2} \leq  B \kappa_i^{- 1}   \Eelinf{p}^{1/p} + B \| \tilde\chi_{ii} - m_{i,\mu_{i}}(D) \tilde\chi_{ii} \|_{L^2},
  \end{align}
  where 
  \begin{align}\label{def:bounditerativeargument2}
B \defas B(g,\tilde\chi_{ii},\lambda) \defas \sup_{a \in \{-3,-1,1,-3 \} - F_{ii}} \sup_{b \in (-3\lambda - F_{ii},3\lambda - F_{ii} )} \vert G(a,b) \vert
\end{align}
and $G$ is a polynomial satisfying the relation  
$	g(a) - g(b) \defas (a-b) G(a,b) $  for $a,b \in \R$.
\end{proposition}
Proposition~\ref{prop:rectanglereduction} shows that for $\mu_i$ significantly larger than $\mu_j$ we have to estimate fewer coefficients with surface energy.
As we have both $g(\chi_{11}) = \chi_{22}$ and $-g(\chi_{22}) = \chi_{11}$ due to our choice of wells, and hence
\begin{align}\label{eq:polyonimalrelationF}
	g(\tilde \chi_{11} + F_{11}) - F_{22} = \tilde \chi_{22} \quad \text{ and } \quad    - g(\tilde \chi_{22}+F_{22}) - F_{11} = \tilde \chi_{11},
\end{align}
this estimate allows for an iteration, successively improving the sharpness of the lower bound, see Section~\ref{sec:proofslower}.
This result particularly improves the results from previous articles as the quantities on the right-hand side appear \emph{linearly}
	instead of with a power strictly smaller than $1$. Moreover, we derived a bound for $B$ as it will be strongly related with the lower scaling bound. 
	Eventually, we note that the factor $(1- \kappa_i^2)^{-1/2}$ will be negligible in our analysis as $\kappa_i$ will be small parameters converging to zero.

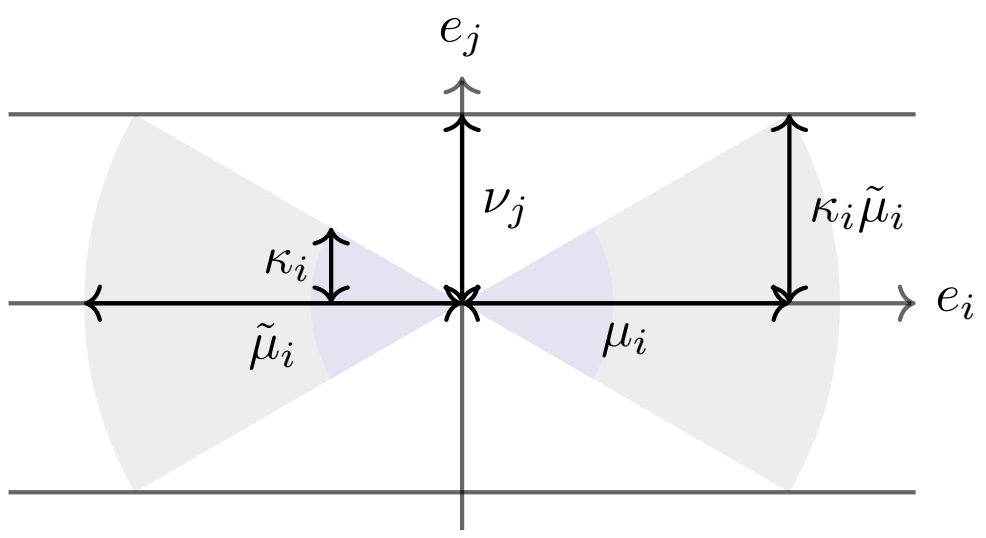
\begin{figure}[t]
	\centering
	\begin{tikzpicture}[thick, scale= 0.5]

\def\mutikz{5}     
\def\kappatikz{30}

\fill[gray!20, opacity=0.7]
  (0,0) --
  ({\mutikz*cos(\kappatikz)},{\mutikz*sin(\kappatikz)}) arc[start angle=\kappatikz, end angle=-\kappatikz, radius=\mutikz] --
  (0,0) -- cycle;

  \fill[gray!20, opacity=0.7]
  (0,0) --
  ({\mutikz*cos(\kappatikz+180)},{\mutikz*sin(\kappatikz+180)}) arc[start angle=\kappatikz+180, end angle=-\kappatikz+180, radius=\mutikz] --
  (0,0) -- cycle;

\fill[blue!20, opacity=0.3]
  (0,0) --
  ({2*cos(\kappatikz)},{2*sin(\kappatikz)}) arc[start angle=\kappatikz, end angle=-\kappatikz, radius=2] --
  (0,0) -- cycle;

  \fill[blue!20, opacity=0.3]
  (0,0) --
  ({2*cos(\kappatikz+180)},{2*sin(\kappatikz+180)}) arc[start angle=\kappatikz+180, end angle=-\kappatikz+180, radius=2] --
  (0,0) -- cycle;

\draw[<->] (0,0) -- (0,{5*sin(\kappatikz)})  node[midway, right]{ $ \nu_j $}; 

\draw[<->] ({5*cos(\kappatikz)},0) -- ({5*cos(\kappatikz)},{5*sin(\kappatikz)})  node[midway, right]{ $ \kappa_i \tilde \mu_i $}; 

\draw[<->] (-5,0) -- (0,0)  node[midway, below]{ $ \tilde \mu_i $};

\draw[<->] (0,0) -- ({5*cos(\kappatikz)},0) node[midway, below]{ $ \mu_i $};

\draw[<->] ({2*cos(\kappatikz+180)},0) -- ({2*cos(\kappatikz+180)},{-2*sin(\kappatikz+180)}) node[midway, left]{ $\kappa_i $};

\draw[->,opacity=0.6] (-6,0) -- (6,0) node[anchor=west] { $  $ };
\draw[->,opacity=0.6] (0,-3) -- (0,3) node[anchor=south] { $ $ };
\draw (6,0) node[right]  {$e_i$};
\draw (0,3) node[above]  {$e_j$};
 \draw[-,opacity=0.6] (-6,2.5) -- (6,2.5) node[anchor=west] { $  $ };
 \draw[-,opacity=0.6] (-6,-2.5) -- (6,-2.5) node[anchor=west] { $  $ };
\end{tikzpicture}
\caption{Visualization of the estimate \eqref{eq:relwithcones} in the proof of Proposition~\ref{prop:rectanglereduction}.
The surface in blue depicts the intersection of the unit disc with the cone.}
\label{fig:cone}
\end{figure}

	 \begin{proof}
Let $\tilde m_{j,\mu_j}(D)$ be the multiplier with coefficients $\tilde m_{j,\mu_j} (\xi) = 1 -m_{j,\mu_j} (\xi)$ which is uniformly bounded in $L^2$ with constant $1$ and vanishes for $\xi$ with $\vert \xi_j \vert \leq \mu_j$. 
Hence, the left-hand side of \eqref{eq:estimatestripered} is equal to $\Vert \tilde m_{j,\mu_j}(D) \tilde \chi_{jj} \Vert_{L^2}$.
We claim that for $\nu_j = \lfloor\frac{\mu_j}{2\deg(g)} \rfloor  $   it holds that
\begin{align}\label{eq:fouriersup}
 \tilde m_{j,\mu_j}(D)     g (     \psi_{j,\nu_j}(D)  \tilde\chi_{ii}  )  = 0.
\end{align}  
Indeed, as $g$ is a polynomial with degree $\deg(g)$, the support of $g(      \psi_{j,\nu_j}(D)\tilde{\chi}_{ii})$ can be estimated by a $\deg(g)$-fold Minkowski sum of the support of $      \psi_{j,\nu_j}(D) \tilde{\chi}_{ii}$. 
(The support of $ \psi_{j,\nu_j}(D) \tilde{\chi}_{ii}$ is given by a stripe of length $4   \nu_j  $, whereas
 $\tilde m_{j,\mu_j}$ vanishes on a stripe of length $2 \mu_j$. Hence, \eqref{eq:fouriersup} is satisfied as long as
$2 \mu_j \geq 4\nu_j \deg(g)$.)

As a consequence of \eqref{eq:fouriersup}, $\tilde \chi_{jj} = g(\tilde \chi_{ii})$, and the $L^2$-boundedness of $\tilde m_{j,\mu_j}$  we find that
\begin{align}\label{eq:Lipschitzproblem}
     \|  \tilde m_{j,\mu_j}(D)  \tilde{\chi}_{jj}\|_{L^2} = \|  \tilde m_{j,\mu_j}(D) (      \tilde{\chi}_{jj}  -  g (     \psi_{j,\nu_j}(D) \tilde{\chi}_{ii} )     ) \|_{L^2}  \leq  \|       g(\tilde \chi_{ii}) -  g (     \psi_{j,\nu_j}(D) \tilde{\chi}_{ii}  )   \|_{L^2}.
\end{align}
Further, it holds that
\begin{align}\label{eq:Lipschitzproblem2}
  \|g(\tilde\chi_{ii} ) - g( \psi_{j,\nu_j}(D)  \tilde\chi_{ii})\|_{L^2} \leq \|\tilde\chi_{ii}  -  \psi_{j,\nu_j}(D) \tilde\chi_{ii} \|_{L^{2}} \| G (\tilde\chi_{ii} ,  \psi_{j,\nu_j}(D) \tilde\chi_{ii}) \|_{L^{\infty}}  .
\end{align}
The latter factor is bounded by
\begin{align*}
	\| G (\tilde\chi_{ii} ,  \psi_{j,\nu_j}(D) \tilde\chi_{ii}) \|_{L^{\infty}}  \leq B(g,\tilde\chi_{ii},\lambda),
\end{align*}
defined in \eqref{def:bounditerativeargument2}.
As visualized in Figure~\ref{fig:cone}, we choose $\nu_j = \lfloor\tilde \mu_i \kappa_i\rfloor$ for $\tilde \mu_i > 0$ and $\kappa_i>0$.
These parameters correspond to the largest cone fitting into the stripe of length $2\nu_j$.
$\kappa_i$ represents the angle of the cone, and all the Fourier coefficients outside the cone can be estimated with elastic energy.
More precisely, Lemma~\ref{lem:Fourierchar} implies that
\begin{align}\label{eq:relwithcones}
	  \| \tilde\chi_{ii} - \psi_{j,\nu_j}(D) \tilde\chi_{ii} \|_{L^2} \leq  \kappa_i^{- 1}   \Eelinf{p}^{1/p} + \| \tilde\chi_{ii} - m_{i,  \mu_i}(D) \tilde\chi_{ii} \|_{L^2},
\end{align}
where $ \mu_i = \tilde \mu_i \sqrt{1- \kappa_i^2}$. The last term corresponds to the Fourier coefficients which are not covered by the estimate with elastic energy. These are all the Fourier coefficients with $\vert \xi_i \vert > \mu_i$.
 (If $\nu_j = 0$, the Fourier series corresponding to $\tilde\chi_{ii} - \psi_{j,\nu_j}(D) \tilde\chi_{ii}$ has non-zero coefficients for $\vert\xi_j \vert \geq 1$, and we may choose $ \kappa_i \tilde \mu_i = 1/2$, leading to $\mu_i = \frac{1}{2} \kappa_i^{-1} \sqrt{1 - \kappa_i^2}$.)
 Eventually, we obtain the relation 
 \begin{align*}
	\left\lfloor \frac{\mu_j}{2 \deg(g)} \right\rfloor = \lfloor \mu_i \kappa_i (1-\kappa_i^2)^{-1/2} \rfloor.
 \end{align*}
If the latter relation is satisfied without floor functions it is also satisfied with them, i.e., our choice of parameters is admissible.
In view of \eqref{eq:Lipschitzproblem}, \eqref{eq:Lipschitzproblem2}, and \eqref{eq:relwithcones}, this concludes the proof.
	 \end{proof}
If $g$ was Lipschitz continuous, the estimate in \eqref{eq:Lipschitzproblem2} became trivial.
This estimate is the core of our argument as it requires us to use multipliers with are bounded in $L^\infty$.

\subsection{Proof of Theorem~\ref{mainth}}\label{sec:proofslower}
In this section, we prove Theorem~\ref{mainth}, following the strategy outlined throughout the previous three subsections.
\begin{proof}[Proof of Theorem~\ref{mainth}]       
	Let $F \in K^{\rm qc} \setminus K$. Then, either $\tilde \chi_{11}$ or $\tilde \chi_{22}$ does not attain the value $0$. By symmetry, it suffices to perform the argument for $\tilde \chi_{11} \neq 0$.
	
	\emph{Step 1 (Initialization as a Fourier series)}:
	For $\mu_1 \in (0,1)$ the Fourier coefficients for $m_{1,\mu_1}(D)\tilde \chi_{11} $ vanish for $\vert \xi_1 \vert \neq 0$. For $\xi_1 = 0$ and arbitrary $\xi_2 $, the Fourier coefficients are controlled via Lemma~\ref{lem:Fourierchar}.
	Using Plancherel's formula and the triangle inequality, we thus obtain a constant $c>0$ such that for all $\kappa \in (0,1)$ it holds that
	\begin{align}\label{eq:classicalfirststep}
		0 < c  = \Vert \tilde \chi_{11} \Vert_{L^2}^2  &\leq \vert \mathcal{F}(\tilde \chi_{11})(0) \vert^2 + \sum\limits_{ \xi \in \mathbb{Z}^2 : \vert \xi_1 \vert \geq \mu_1} \vert \mathcal{F}(\tilde \chi_{11})(\xi) \vert^2 + \sum\limits_{ \xi \in \mathbb{Z}^2 : \vert \xi_1 \vert =0 , \ \xi_2 \neq 0 } \vert \mathcal{F}(\tilde \chi_{11})(\xi) \vert^2   \notag \\
		  &\leq \Vert \tilde \chi_{11} - m_{1,\mu_1}(D) \tilde \chi_{11} \Vert_{L^2}^2 + \kappa^{-2} \Eelinf{p}^{2/p}.
	\end{align}

\emph{Step 2 (The bootstrap argument)}:
As the right-hand side of \eqref{eq:classicalfirststep} features the term $\Vert \tilde \chi_{11} - m_{1,\mu_1}(D) \tilde \chi_{11} \Vert_{L^2}$, our next goal is to apply Proposition~\ref{prop:rectanglereduction} and to iterate the latter.

\emph{Step 2.1 (The polynomial)}:
It holds that $g(\chi_{11}) = \chi_{22}$ for
$g(t) = - \frac{5}{12} t^3 + \frac{41}{12} t$, which is a polynomial of degree $3$. Moreover, $g(a)- g(b) = (a-b) G(a,b)$ for $G(a,b) = \frac{41}{12}- \frac{5}{12}(a^2 + ab +b^2)$.
Passing from $\chi$ to $\tilde \chi$, see \eqref{eq:polyonimalrelationF}, we hence observe
$\tilde g(\tilde \chi_{11}) \defas  g(\tilde \chi_{11} + F_{11}) - F_{22} = \tilde \chi_{22}$
as well as
$$\tilde g(a) - \tilde g(b) = g(a+ F_{11}) - g(b + F_{11}) = (a-b) G(a+ F_{11}, b + F_{11}).$$
Let us determine the value of $B$ in \eqref{def:bounditerativeargument2} and compute the maximum of $\vert G(a+ F_{11}, b + F_{11}) \vert$.
As $a^2 + ab + b^2 \geq 0$ for all $a,b \in \R$, the maximum of $G$ is attained at $(0,0)$ with value $\frac{41}{12}$.
It remains to compute the minimum. Recall that
$\lambda$ is defined in \eqref{eq:bestbound}.
According to \eqref{def:bounditerativeargument2}, the admissible values for $a$ are $-3- F_{11}, -1- F_{11},1- F_{11},3-F_{11}$. The range of $b$ is $ [-F_{11}- 3\lambda,-F_{11}+ 3 \lambda ]$.
The maximal value of $\vert G(a+ F_{11}, b + F_{11}) \vert$ hence does not depend on $F_{11}$, and we may set $F_{11} = 0$.
The extremal point of $a^2 + ab + b^2$ is attained for $ a = 3$ and $b = 3\lambda$.
Hence, we have $B = \max\{\frac{41}{12},  \frac{5}{12}(9(1+\lambda + \lambda^2)) - \frac{41}{12} \} $.
Recall that the value of $\lambda$ was defined in \eqref{eq:bestbound}, and we have
\begin{align*}
	B = 2 + \frac{25 \sqrt{3}}{2\pi} + \frac{45}{\pi^2} \approx 13.45 .
\end{align*}
\emph{Step 2.2 (The estimate)}:
By an iteration of Proposition~\ref{prop:rectanglereduction}, we hence derive that
	 \begin{align*}
		\Vert  \tilde \chi_{11} - m_{1,\mu_1}(D) \tilde \chi_{11} \Vert_{L^2}  & \leq \sum_{k=0}^{m-1}  B^k  \kappa^{- 1}   \Eelinf{p}^{1/p} + B^m \Vert  \tilde \chi_{jj} - m_{j,\mu_{m}}(D) \tilde \chi_{jj} \Vert_{L^2} 
	 \end{align*}
	 for some $j \in \{1,2\}$ and 
	 \begin{align}\label{eq:relationiteration}
	 \mu_1 = 6^m  \kappa^m (1- \kappa^2)^{-m/2} \mu_m
	 \end{align}
	 as $2 \deg(g) = 6$.
The estimate $\sum_{k=0}^{m-1}  B^k \leq (B-1)^{-1} B^m$ and Lemma~\ref{lem:highfrequencyoutsidecones} imply that
\begin{align}\label{iteration2}
	\Vert  \tilde \chi_{11} - m_{1,\mu_1}(D) \tilde \chi_{11} \Vert_{L^2}  &\leq C B^m \kappa^{- 1 }   \Eelinf{p}^{1/p}  \notag \\
	& \qquad + C B^m \big( \mu_{m}^{-1}  (\varepsilon^{-1} \varepsilon \Esurf + \mu_m^{-1} \Per ) \big)^{1/2} .
\end{align}

\emph{Step 3 (Optimization)}:
We combine the estimates \eqref{eq:classicalfirststep} and \eqref{iteration2} and use Young's inequality (in the form $ab \leq \delta^{-1}q^{-1}a^{q} + \delta \tilde q^{-1} b^{\tilde q})$ such that we obtain a smaller constant $c>0$ satisfying
\begin{align*} 
	c (1- B^{2m} \mu_m^{-1} )\leq \left(B^{mp}  \kappa^{- p }   +  B^m    \mu_{m}^{-2}  \varepsilon^{-2} \right) \Energyinf{p}   .
\end{align*}
As $\mu_1 \in (0,1)$ fixed, the relation in \eqref{eq:relationiteration} implies that
\begin{align*}
	B^{2m} \mu_m^{-1} \leq 2 B^{2m} 6^m \kappa^m < \frac{1}{2}
\end{align*}
for $\kappa$ sufficiently small, and hence
\begin{align*}
	c \leq \left( M^m  \kappa^{- p }   +  M^m \kappa^{2m}  \varepsilon^{-2} \right) \Energyinf{p}   
\end{align*}
for $M=\max\{B  6^{2},B^{p} \}$. The choice
$\kappa = \varepsilon^{\frac{2}{p+2m}}$ gives 
\begin{align*}
	c \leq  M^m  \varepsilon^{\frac{-2p}{p+2m}} \Energyinf{p} = \exp \big( m \log( M ) + \frac{2p}{p+2m}\log( \varepsilon^{-1}  ) \big)  \Energyinf{p}   .
\end{align*}
We minimize the exponent in $m$ and obtain the choice
$m = - \frac{1}{4} + \sqrt{\frac{1}{16} + p \log( \varepsilon^{-1}  ) \log( M )^{-1}}$ as a minimizer.
As $m$ needs to be a natural number, we therefore select $m^* = \lceil m \rceil$. 
Due to the estimate
\begin{align}\label{eq:assistchat}
	p +2m^* &= p - 1/2 + 2 \sqrt{\frac{1}{16} + p \log( \varepsilon^{-1}  ) \log( M )^{-1}} \geq p-1/2+ 2\sqrt{p \log( \varepsilon^{-1}  ) \log( M )^{-1}} \notag \\
	& \geq 2 \sqrt{p \log( \varepsilon^{-1}  ) \log( M )^{-1}}
\end{align}
we find that
\begin{align*}
	\frac{2p}{p+2m^*}\log( \varepsilon^{-1}  )   \leq \sqrt{p} \log(M)^{1/2}\log( \varepsilon^{-1}  )^{1/2}.
\end{align*}
Thus, as $m^* \leq m+1$ we obtain
\begin{align*}
	c \leq  \exp \big( 2\sqrt{p} \log(M)^{1/2}\log( \varepsilon^{-1}  )^{1/2}\big)  \Energyinf{p} .
\end{align*}
An explicit computation shows
\begin{align*}
2\sqrt{p} \log(M)^{1/2} &  = \max\{ 2 \sqrt{ \log\left(72 + \frac{450 \sqrt{3}}{\pi}+ \frac{1620}{\pi^2} \right) } \sqrt{p}, 2 \sqrt{ \log\left( 2+ \frac{25 \sqrt{3}}{2\pi}+ \frac{45}{\pi^2} \right) } p \} \\
& \leq \max\{ 4.9730 \sqrt{p}, 3.2243 p\}.
\end{align*}
This concludes the proof.
\end{proof}

\begin{remark}[Comments on the sharpness of $C_{\rm low}$]\label{rem:optimalconst}
	The proof shows that the constant heavily depends on $B$ from Proposition~\ref{prop:rectanglereduction}.
	On the one hand, the latter depends on $\lambda$, which is defined as the $L^1$-norm of the kernel for the cut-off of Fourier coefficients.
In principle, a multiplier with better bounds may lead to a smaller constant $C_{\rm low}$. However, the length of the cut-off layer
also influences the constant $C_{\rm low}$:
for example, if one interpolates in \eqref{cut-offmultiplier} over a longer distance, i.e., if $\hat V_N^k = \hat V_N$ for $\vert \xi_1 \vert \leq N $ and $\hat V_N^k = 0$ for $\vert \xi_1 \vert \geq kN $, $k \in \{2, 3, \dots\}$ with a linear interpolation in between, one has $\Vert V_N^k \Vert_{L^1} \leq \frac{k+1 }{k-1}$, see \cite{MEHTA15}.
	This shows that the constant approaches $1$ for large $k$, but the length of the cut-off converges to $\infty$.
	In order to have disjoint Fourier supports in \eqref{eq:fouriersup}, the relation between $\mu_i$ and $\mu_j$ in Proposition~\ref{prop:rectanglereduction} needs to be adjusted, leading to $\mu_j = k \deg(g) \mu_i \kappa_i$. This, however, enforces a loss in the constant $C_{\rm low}$ as it would contribute via $M = \max\{ B k^2 \deg(g)^2, B^p\}$ in the proof of Theorem~\ref{mainth}.
Thus, one seeks to find multipliers with a short truncation range and $L^\infty$-bounds which are close to $1$.
On the other hand, the constant strongly depends on the degree of the polynomial. A higher degree influences $M$ quadratically, but at the same time the constant $B$ might be improved. 
Note, however, that the optimal choice $M$ cannot exceed $\sqrt{2p\log(2)} = 2 \sqrt{p} \log(M)^{1/2}$ due to the upper bound in Theorem~\ref{mainth:upper}.
\end{remark}

\section{Upper bound}
\label{sec:upper} 
The goal of this section is to provide an iterated branching construction, showing that
the constant $C_{\rm upp}$ in Theorem~\ref{mainth:upper} can be chosen within the range $C_{\rm upp} \in (0,(2p)^{1/2} \log^{1/2}(2))$.
We rely on the construction in \cite[Appendix~A]{RT23}, where nested branching constructions for Dirichlet data of finite order have been derived.
In our context, we also have to optimize the number of iterations.
This optimization is inspired by \cite{W97}, where an iterated laminate construction for $p = 2$ was constructed, leading to the range $C_{\rm upp} \in (0,\sqrt{2} \log^{1/2}(2))$ for $p = 2$.
Once the explicit construction is provided in Section~\ref{sec:explicitconst}, the optimal length scales are discussed in Section~\ref{sec:updet}, proving Theorem~\ref{mainth:upper}.

\subsection{Explicit construction}\label{sec:explicitconst}
We divide the construction into three parts. First, we collect the results from \cite[Subsection~3.1]{RT23} and present them in Subsection~\ref{sec:preparationupperbound} in a form which is convenient for the construction of the competitor.
Then, in Subsection~\ref{sec:initalization}, we take the structure of Tartars square into account and initialize the construction by oscillating between the matrices $A_i$ and $J_i$ for some $i \in \{1,\dots, 4\}$.
As much elastic energy is paid when the gradient is close to $J_i$, we replace these parts with further branching constructions between wells $A_{i+1}$ and $J_{i+1}$, successively reducing the elastic energy.
This iteration is provided in Subsection~\ref{sec:theiteration}. 

\subsubsection{Preparations}\label{sec:preparationupperbound}
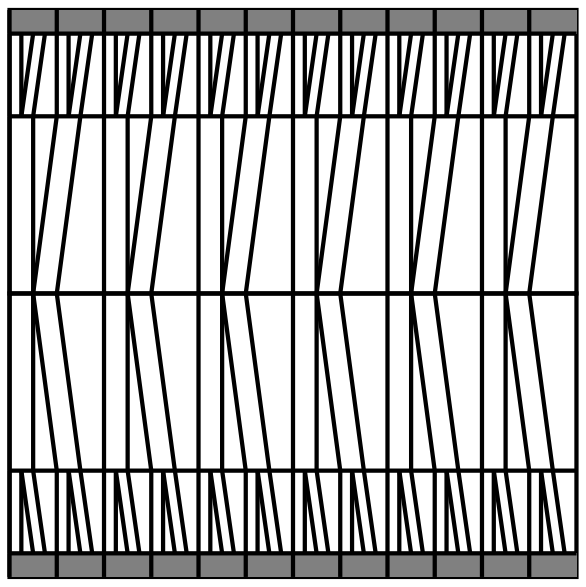
\begin{figure}[t]
\centering 
\begin{tikzpicture}[scale=0.3125]
  \begin{scope} 
\draw[thick] (0,0) rectangle (12,6);
\fill[gray] (0,5.5) rectangle (12,6);
\draw[thick] (0,5.5) -- (12,5.5);  
\draw[thick] (0,3.75) -- (12,3.75);  
\foreach \i in {0,2,...,10} {
    \draw[thick]  (\i, 3.75) -- (\i, 0);
}
\foreach \i in {0,1,...,11} {
    \draw[thick]  (\i, 5.5) -- (\i, 3.75);
}
\foreach \i in {0,1,...,11} {
    \draw[thick]  (\i, 6) -- (\i, 5.5);
}
\foreach \i in {0,2,...,10} {
   \draw[thick] (\i+0.5,0) -- (\i+0.5,3.75);  
   \draw[thick] (\i+0.5,0) -- (\i+1,3.75);  
    \draw[thick] (\i+1,0) -- (\i+1.5,3.75);   
 
}
\foreach \i in {0,1,...,11} {
    \draw[thick] (\i+0.25,3.75) -- (\i+0.25,5.5);  
   \draw[thick] (\i+0.25,3.75) -- (\i+0.5,5.5);  
    \draw[thick] (\i+0.5,3.75) -- (\i+0.75,5.5);  
}
\end{scope}
\begin{scope}[yscale=-1]  
\draw[thick] (0,0) rectangle (12,6);
\fill[gray] (0,5.5) rectangle (12,6);
\draw[thick] (0,5.5) -- (12,5.5);  
\draw[thick] (0,3.75) -- (12,3.75);  
\foreach \i in {0,2,...,10} {
    \draw[thick]  (\i, 3.75) -- (\i, 0);
}
\foreach \i in {0,1,...,11} {
    \draw[thick]  (\i, 5.5) -- (\i, 3.75);
}
\foreach \i in {0,1,...,11} {
    \draw[thick]  (\i, 6) -- (\i, 5.5);
}
 \foreach \i in {0,2,...,10} {
   \draw[thick] (\i+0.5,0) -- (\i+0.5,3.75); 
   \draw[thick] (\i+0.5,0) -- (\i+1,3.75); 
   \draw[thick] (\i+1,0) -- (\i+1.5,3.75); 
}
\foreach \i in {0,1,...,11} {
    \draw[thick] (\i+0.25,3.75) -- (\i+0.25,5.5); 
    \draw[thick] (\i+0.25,3.75) -- (\i+0.5,5.5); 
    \draw[thick] (\i+0.5,3.75) -- (\i+0.75,5.5); 
}
\end{scope}
\end{tikzpicture}
  \caption{ Visualization of a first-order construction.  The picture shows two layers of building blocks and the gray part corresponds to the cut-off layer.}
    \label{fig:selfsimilarconst}
\end{figure}

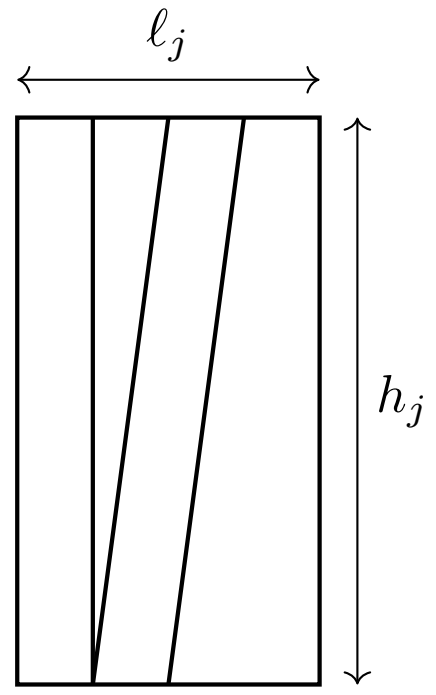
\begin{figure}[t]
\centering 
\hspace{0.03\textwidth}

 \begin{tikzpicture}[scale=1]
  \begin{scope}  
\draw[thick] (0,0) rectangle (2,3.75); 
\foreach \i in {0} {
   \draw[thick] (\i+0.5,0) -- (\i+0.5,3.75); 
   \draw[thick] (\i+0.5,0) -- (\i+1,3.75); 
   \draw[thick] (\i+1,0) -- (\i+1.5,3.75); 
}

\draw[<->] (2.25,0) -- (2.25,3.75) node[midway,right] {$h_j$};  
\draw[<->] (0,4) -- (2,4) node[midway,above] {$\ell_j$};

\end{scope}

\end{tikzpicture}
  \caption{Visualization of a building block. $\partial_1 u_1 $ oscillates among two values such that $u_1$ is piecewise affine and continuous on the four subsets $\omega^{j,k}_1, \dots, \omega^{j,k}_4$.}
    \label{fig:branchingpicture}
\end{figure}

We first recall the construction of a first-order branching construction oscillating between two rank-one connected wells, as it was carried out in \cite[Subsection~3.1]{RT23}.
Let $p\in[1,\infty)$, $0<L, H\le1$, $R=[0,L]\times[0,H]$ and $N\in\N$ such that $N>\frac{4L}{H}$.
Let $A,B\in \R^{2 \times 2}_{\rm diag}$ be such that $A-B =|A-B|e_1\otimes e_1$, $F=\frac{1}{2} ( A+B)$. (The case $F=\lambda A+(1-\lambda)B$ for some $0<\lambda<1$ can also be considered, we stick to the case $\lambda = 1/2$ to simplify the notation.)
The competitor is constructed in a self-similar way, see Figure~\ref{fig:selfsimilarconst}.
As visualized in Figure~\ref{fig:branchingpicture}, we define the width and height of a building block by
\begin{equation}\label{lengthscales}
\ell_j=\frac{L}{2^j N} \quad  \text{ and } \quad 
h_j=\theta^j\frac{H(1-\theta)}{2} \quad \text{ for } \quad 
\theta\in\begin{cases}
(2^{-\frac{p}{p-1}},\frac{1}{2}) & p>1, \\
(0,\frac{1}{2}) & p=1.
\end{cases}
\end{equation}
More precisely, $R$ is divided into building blocks $\omega^{j,k}$ which are defined by
\begin{equation}\label{eq:cells-bra}
\omega^{j,k}:=\begin{cases}\displaystyle
(k\ell_j,y_j)+[0,\ell_j]\times[0,h_j] & j=0,\dots,j^{\rm cut},\quad k=0,\dots, N2^j-1, \\ 
\\ \displaystyle
[0,L]\times\Big[y_{j^{\rm cut}+1},y_{j^{\rm cut}+1}+\frac{H}{2}\theta^{j^{\rm cut}+1}\Big] & j=j^{\rm cut}+1,
\end{cases}
\end{equation} 
where $\omega^{j^{\rm cut}+1,k} = \omega^{j^{\rm cut}+1} $ does not depend on $k$.
Here, we set $y_{j}=H-\frac{H}{2}\theta^{j} $ such that  $h_{j}=y_{j+1}-y_{j}=\theta^{j}\frac{H(1-\theta)}{2}$ for $j = 0, \dots, j^{\rm cut}$.
Further, $y_0 = \frac{H}{2}$, i.e., the building blocks in \eqref{eq:cells-bra} only fill the rectangle $[0,L] \times [\frac{H}{2},H]$ (up to the cut-off layer),
and the lower rectangle is derived by reflection, see Figure~\ref{fig:selfsimilarconst}. 
The choice of $\theta$ will ensure that the energy among every building block is summable, see \eqref{boundEwell} and \eqref{boundEsurf} below.
As $\theta < 1/2$, there is a last index for which $\ell_{j_1} \leq h_{j_1}$, denoted by $j^{\rm cut}_{1} $.
This condition defines the cut-off region $[0,L]\times\big[y_{j^{\rm cut}+1},H\big]$.
 To simplify the notation, the height of the cut-off layer $h_{j^{\rm cut}+1}$ is defined as $h_{j^{\rm cut}+1}= \frac{\theta}{1-\theta} h_{j^{\rm cut}}$. This is not consistent with the formula in \eqref{lengthscales}, however, $h_{j^{\rm cut}+1} \sim \theta^{j^{\rm cut}+1}\frac{H(1-\theta)}{2}$, i.e., the quantities are comparable, which suffices for the computation below.
For $j = 1, \dots, j^{\rm cut}$, the building blocks are subdivided into the sets
\begin{align*}
\omega^{j,k}_1 &= (k\ell_j,y_j)+\Big\{(x_1,x_2)\in [0,\ell_j]\times[0,h_j] \,:\, 0<x_1<\frac{\ell_j}{4}\Big\},\\
\omega^{j,k}_2 &= (k\ell_j,y_j)+\Big\{(x_1,x_2)\in [0,\ell_j]\times[0,h_j] \,:\, \frac{\ell_j}{4}<x_1<\frac{\ell_j}{4}+\frac{\ell_j\,x_2}{4h_j}\Big\},\\
 \omega^{j,k}_3&  = (k\ell_j,y_j)+\Big\{(x_1,x_2)\in [0,\ell_j]\times[0,h_j] \,:\, \frac{\ell_j}{4}+\frac{\ell_j\,x_2}{4h_j}<x_1<\frac{1}{2}\ell_j+\frac{\ell_j\,x_2}{4h_j}\Big\},\\
\omega^{j,k}_4 &= (k\ell_j,y_j)+\Big\{(x_1,x_2)\in [0,\ell_j]\times[0,h_j] \,:\, \frac{1}{2}\ell_j+\frac{\ell_j\,x_2}{4h_j}<x_1<\ell_j\Big\}.
\end{align*}
This is highlighted in Figure~\ref{fig:branchingpicture}.
 Let us define
 the surface energy on one building block as $E_{\rm surf} (v;\omega^{j,k}) = \Vert D \, {\rm Pr}_K (\nabla v) \Vert_{TV(\omega^{j,k})}$, where ${\rm Pr}_K$ denotes the projection of $\nabla v$ onto $K$.
 Then, the results from \cite{RT23} show the following.
\begin{lemma}[Energy contribution of a first-order construction]\label{lem:firstorderconstr}
There exists a function $v\in W^{1,\infty}(R;\R^2)$ with the following properties:
\begin{itemize}
	\item[(i)] For $j\in\{0,\dots,j^{\rm cut}\}$ we have
\begin{equation}\label{eq:elbuildingblock}
 \int_{\omega^{j,k}} \dist^p(\nabla v, \{A, B\}) \di x  \lesssim \frac{\ell_j^{p+1}}{h_j^{p-1}} \qquad \text{and} \qquad  E_{\rm surf} (v;\omega^{j,k}) \lesssim h_j 
\end{equation}
for every  $k\in\{0,\dots,k_j-1\}$ with $k_j:=N2^j$.
On $\omega^{j,k}_2$ and $\omega^{j,k}_4$ we have ${\rm Pr}_K (\nabla v) = B$, and on $\omega^{j,k}_1$ and $\omega^{j,k}_3$ we have ${\rm Pr}_K (\nabla v) = A$.
\item[(ii)] For $j = j^{\rm cut}+1$ it holds that
\begin{align}\label{eq:cut-off}
\int_{\omega^{j^{\rm cut}+1}} \dist^p(\nabla v, \{A, B\}) \di x &\lesssim  h_{j^{\rm cut}+1}^{1-p} \ell_{j^{\rm cut} +1}^{p+1} N 2^{j^{\rm cut}+1}  ,\notag\\
   E_{\rm surf} (v;\omega^{j^{\rm cut}+1}) &\lesssim h_{j^{\rm cut}+1} N 2^{j^{\rm cut}+1}.
\end{align}
\item[(iii)] $v(x) = Fx$ on $\partial R$.
\end{itemize}
\end{lemma}
\begin{proof} 
	The result is proved in
	\cite[Subsection~3.1, Lemma~3.1 and Lemma~3.2]{RT23}. As we employ the same notation, we omit the proof. The function $v\in W^{1,\infty}(R;\R^2)$ is constructed in a self-similar way such that it is piecewise affine on $\omega^{j,k}$ 
for every $j\in\{0,\dots,j^{\rm cut}\}$, $k\in\{0,\dots,k_j-1\}$ with $k_j:=N2^j$ if $j\le j^{\rm cut}$. 
	Regarding the cut-off layer, we note that gradient on the cut-off layer is uniformly bounded. Hence, the elastic energy is bounded by $h_{j^{\rm cut}+1} \ell_{j^{\rm cut} +1} N 2^{j^{\rm cut}+1} \leq h_{j^{\rm cut}+1}^{1-p} \ell_{j^{\rm cut} +1}^{p+1} N 2^{j^{\rm cut}+1} $ as $\ell_{j^{\rm cut} +1} > h_{j^{\rm cut} +1}$, see \cite[Proof of Lemma 3.2]{RT23}.
	In a similar spirit, the surface energy is controlled by $h_{j^{\rm cut}+1} N 2^{j^{\rm cut}+1}$.
\end{proof}
	If one consider matrices $A,B\in \R^{2 \times 2}_{\rm diag}$ such that $A-B =|A-B|e_2\otimes e_2$ and $F=\frac{1}{2} (A+B)$.
Then, the roles of length and height in the construction above are reversed, and one needs to branch into the perpendicular direction.
 
\subsubsection{Initialization}\label{sec:initalization}
Let us turn to the proof of Theorem~\ref{mainth:upper}. We carry out our construction for $\Omega = (0,1)^2$.
We define the auxiliary matrices
				\begin{align}\label{eq:auxmatrices}
		  J_1 = \diag(-1,1) , \ J_2 = \diag(1,1) , \ J_3 = \diag(1,-1), \ J_4 = \diag (-1,-1), 
	\end{align}
	and observe that the associated the rank-one connections with the wells are given as
	\begin{align*} 
		J_1 - A_1 = A_3 - J_3 = 4 e_2 \otimes e_2, \qquad  J_2 - A_2 =  A_4 - J_4 = 4 e_1 \otimes e_1,
	\end{align*}
see also Figure~\ref{fig:Tartarsquare}.
Furthermore, we provide the construction for the Dirichlet boundary conditions with $F = \frac{1}{2} (A_2 + J_2)$. (This choice is convenient as in every step we will be split the parts where we pay lots of elastic energy into two parts with equal volume fraction, but we note that a general $F \in K^{\rm qc}$ can also be handled with minor adaptations.)
Further, we introduce the length scales $r_1>r_2> \dots$ with $\lim_{\varepsilon \to 0} r_i r_{i-1}^{-1} = 0$ for all $i\in \N$.
In view of Lemma~\ref{lem:firstorderconstr}, we consider a branching construction between the rank-one connected wells $A_2$ and $J_2$, and set
\begin{align}\label{eq:lengthdef}
	\ell_{j_1} = \frac{1}{2^{j_1} N_1}, \qquad 
h_{j_1} =\theta^{j_1}\frac{(1-\theta)}{2}
\end{align}
with
\begin{align*}
	j_1=0,\dots,j^{\rm cut}_{1} ,\quad k_1=0,\dots, N_1 2^{j_1}-1,
\end{align*}
where we choose $N_1  \defas \lceil r_1^{-1} \rceil$.

\begin{remark}[Current energy contribution]\label{rem:firstorderenergycont}
For convenience, we introduce some notation which will be convenient for the iteration below, and compute the current energy contribution.
The energy contribution consists of four parts.
\begin{itemize}
	\item[(1)] In the half of the building blocks few elastic energy is paid, defined by $E_{\rm well}^{(j_1,k_1)}$. More precisely, on $\omega_{2}^{j_1,k_1}$ and $\omega_{4}^{j_1,k_1}$, we have ${\rm Pr}_K (\nabla v) = A_2$, and \eqref{eq:elbuildingblock} gives
	\begin{align*}
		E_{\rm well}^{(j_1,k_1)} \lesssim \frac{\ell_{j_1}^{p+1}}{h_{j_1}^{p-1}}
	\end{align*}
	on every building block $\omega^{j_1,k_1}$.
	\item[(2)] In the other parts of the building blocks much elastic energy is paid, denoted by $E_{\rm pay}^{(j_1,k_1)}$, as ${\rm Pr}_K (\nabla v) = J_2$. This part scales like the volume of the phase as $J_2 \notin K$. More precisely, we have
	\begin{align*}
		E_{\rm pay}^{(j_1,k_1)} \lesssim \left( \vert \omega_{2}^{j_1,k_1} \vert + \vert \omega_{4}^{j_1,k_1} \vert \right) \dist^p(J_2, K).
	\end{align*}
	\item[(3)] The surface energy on the building blocks is denoted by $E_{\rm surf}^{(j_1,k_1)}$, which is controlled by $h_{j_1}$, see \eqref{eq:elbuildingblock}.
	\item[(4)] Eventually, $E^{(j_1^{\rm cut})}_{\rm cut}$ denotes the energy on the whole cut-off part, and is estimated via \eqref{eq:cut-off}.
\end{itemize}
Denoting by $E^{(1\rm st)}$ the total energy contribution of the first-order branching construction, we find that
\begin{align*}
	E^{(1\rm st)} &= 2 E^{(j_1^{\rm cut})}_{\rm cut} + 2 \sum_{j_1=0}^{j^{\rm cut}_{1}} \sum_{k_1=0}^{N_1 2^{j_1}-1}  \left(E_{\rm well}^{(j_1,k_1)} + E_{\rm pay}^{(j_1,k_1)} + \varepsilon E_{\rm surf}^{(j_1,k_1)}  \right) \\
	&\lesssim  \sum_{j_1=0}^{j^{\rm cut}_{1}+1} N_1 2^{j_1}  \left(     \frac{\ell_{j_1}^{p+1}}{h_{j_1}^{p-1}} +  \varepsilon h_{j_1} \right) +  \sum_{j_1=0}^{j^{\rm cut}_{1}+1} \sum_{k_1=0}^{N_1 2^{j_1}-1}  \ell_{j_1} h_{j_1}  \\
		&\lesssim   \sum_{j_1=0}^{j^{\rm cut}_{1}+1} N_1 2^{j_1}  \left(  \frac{\ell_{j_1}^{p+1}}{h_{j_1}^{p-1}} +  \varepsilon h_{j_1} \right) +  \frac{1}{2}  , 
\end{align*}
where the cut-off part $E^{(j_1^{\rm cut})}$ is a lower order term and controlled by the series due to \eqref{eq:cut-off}.
The factor $2$ in the estimate above arises as our notation only covers the upper part $[0,1] \times [1/2,1]\subset \Omega$, and the lower part is derived by reflection, see Figure~\ref{fig:selfsimilarconst}.
 \end{remark}
\subsubsection{The iteration}\label{sec:theiteration}
\begin{figure}[t]
    \centering
    \includegraphics[width=\textwidth]{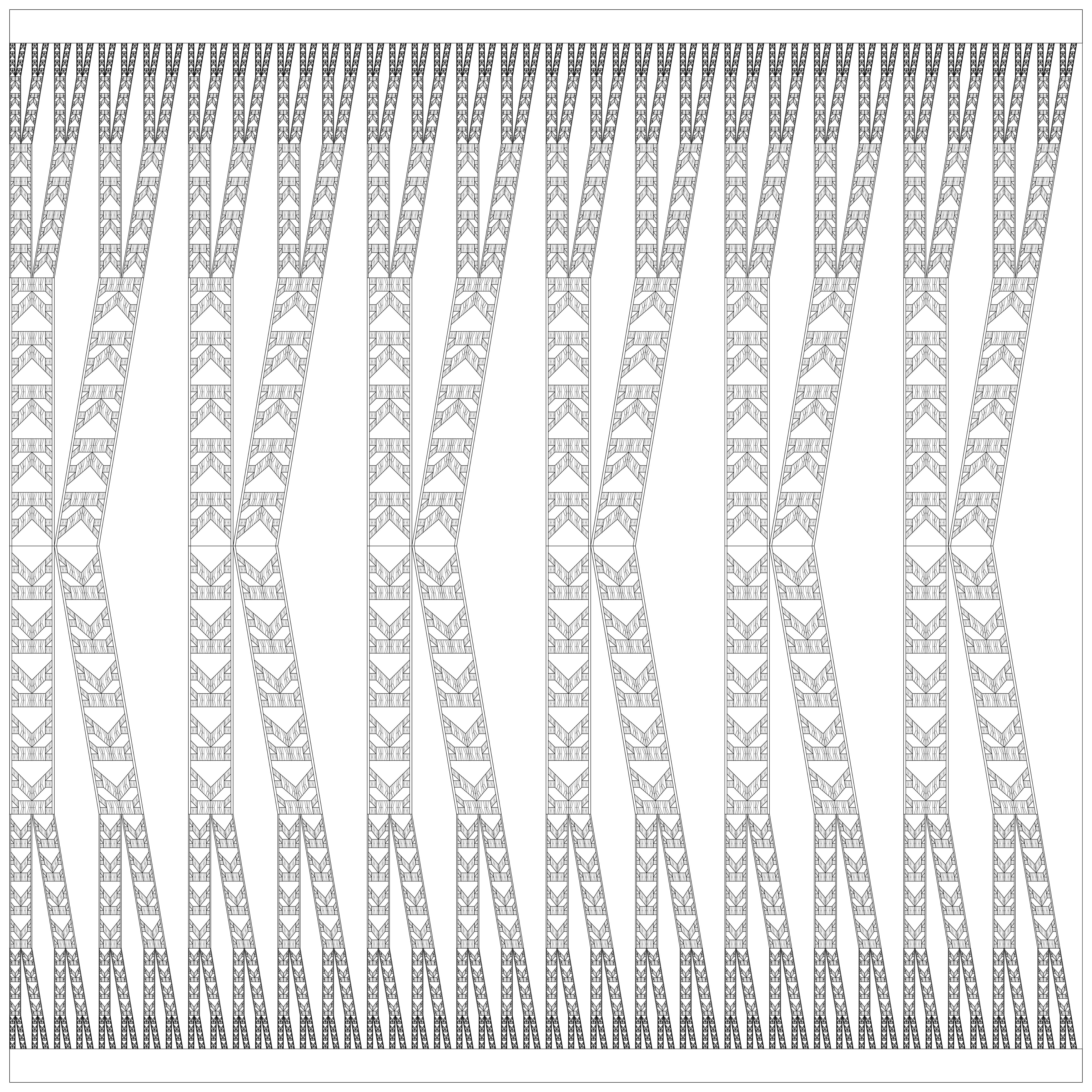}
    \caption{Visualization of an iterated branching construction. The picture schematically displays a branching construction of third order, whereas the order increases in our construction as $\varepsilon \to 0$.}
    \label{fig:branching}
\end{figure}

    As visualized in Figure~\ref{fig:branching}, we perform iterated branching constructions.
\subsubsection*{Second-order construction}
Since we pay elastic energy, where ${\rm Pr}_K(\nabla v) = J_2$, we add in these parts a further branching construction oscillating between $A_3$ and $J_3$, recalling Figure~\ref{fig:Tartarsquare}.
This time, the compatibility direction is $e_2$ instead of $e_1$. Hence, we perform a branching construction with oscillations in the $e_2$-direction.
More specifically, for every $j_1=0,\dots,j^{\rm cut}_{1} $ and $ k_1=0,\dots, N_1 2^{j_1}-1$ the added constructions in the blocks $\omega^{j_1,k_1}_1$ and $\omega_3^{j_1,k_1}$ 
generate the blocks
$\omega^{j_1,j_2,k_1,k_2}_1$ and $\omega_3^{j_1,j_2,k_1,k_2}$. (Note that this is also carried out in $[0,1] \times [0,1/2]$ by reflection.)
The width and height of these blocks are given as
\begin{align*}
	\ell_{j_2} = \frac{h_{j_1}}{2^{j_2} N_2}, \qquad  
h_{j_2} =\theta^{j_2}\frac{\ell_{j_1}(1-\theta)}{8},
\end{align*}
with
\begin{align*}
	j_2=0,\dots,j^{\rm cut}_{2} ,\quad k_2=0,\dots, N_2(j_1) 2^{j_2}-1.
\end{align*}
 $\ell_{j_2}$ describes the width of a building block, and the associated edge again corresponds to the direction of oscillations.
 At this point, we stress that $\ell_{j_2}$ and $h_{j_2}$ depend on the previous length scales $\ell_{j_1}$ and $h_{j_1}$. More specifically, in what follows, the introduced length scales will always depend on the length scales in the previous layers. For convenience of the reader, we omit this notation.
 The factor $8$ arises as the width of $\omega^{j_1,j_2,k_1,k_2}_1$ and $\omega_3^{j_1,j_2,k_1,k_2}$ is equal to $\ell_{j_1}/4$. Hence, the height $h_{j_2}$ is given as $\frac{1}{4} \theta^{j_2}\frac{\ell_{j_1}(1-\theta)}{2}$.
Here, we note that the branching constructions in Lemma~\ref{lem:firstorderconstr} are carried out on a rectangle. We will use that the energy contribution for the replacements in $\omega_3^{j_1,k_1}$ is comparable to the contribution in $\omega_1^{j_1,k_1}$ thanks to a simple shear and a change of variables, see \cite[Equation (42)]{RT23}.

\subsubsection*{Third-order construction}
As we still pay too much elastic energy in the case ${\rm Pr}_K(\nabla v) = J_3$, we replace these parts with branching constructions oscillating between $A_4$ and $J_4$ into the direction $e_1$.
Again, we define new blocks  $\omega^{j_1,\dots,j_3,k_1,\dots, k_3}$ and $\omega^{j_1,\dots,j_3,k_1,\dots, k_3}$ with width and height
\begin{align*}
	\ell_{j_3} = \frac{h_{{j_{2}}}}{2^{j_3} N_3}, \qquad   
h_{{j_3}}=\theta^{{j_3}}\frac{\ell_{{j_{2}}}(1-\theta)}{8}
\end{align*}
with
\begin{align*}
	{j_3} =0,\dots,j^{\rm cut}_{3} ,\quad k_3=0,\dots, N_3 2^{{j_3}}-1.
\end{align*}

 \subsubsection*{$m$th-order construction}
 In this way, we exploit the $T_4$-structure by iterating this algorithm, see Figure~\ref{fig:Tartarsquare}.
In the parts, where ${\rm Pr}_K(\nabla v) = J_j$ for some $j = 1,\dots, 4$, we insert further branching constructions reducing the elastic energy.
Given $\ell_{j_{m-1}  }$, $h_{j_{m-1}  }$ and blocks   $\omega^{j_1,\dots,j_{m-1},k_1,\dots, k_{m-1}}_1$ and $\omega^{j_1,\dots,j_{m-1},k_1,\dots, k_{m-1}}_3$, we define new blocks  $\omega^{j_1,\dots,j_m,k_1,\dots, k_m}$ and $\omega^{j_1,\dots,j_m,k_1,\dots, k_m}$ with
 width and height
\begin{align}\label{def:lengthsblocksit}
	&\ell_{j_m} = \frac{h_{{j_{m-1}}}}{2^{j_m} N_m}, \quad h_{{j_m }}=\theta^{{j_m }}\frac{\ell_{{j_{m-1} }}(1-\theta)}{8},
\end{align}
where
\begin{align*}
	{j_m} =0,\dots,j^{\rm cut}_{m} ,\quad k_m=0,\dots, N_m 2^{{j_m}}-1.
\end{align*}
Recall that the width and height $\ell_{j_m}$ and $h_{j_m}$ also depend on the indices $j_1, \dots, j_{m-1}$.
Eventually, we find the following lemma.
\begin{lemma}[$m$th-order construction]\label{lem:preliminaryiteration}
	The $m$th iteration of the previous algorithm generates an energy contribution of the form
	\begin{align}\label{energycontribution}
		 E^{(m\rm th)} & = \sum_{i = 0}^{m-1}  \left(   \sum_{j_0=0}^{j^{\rm cut}_{0}} \sum_{k_0=0}^{N_0 2^{j_0}-1}   \dots     \sum_{j_i=0}^{j^{\rm cut}_{i}} \sum_{k_i=0}^{N_i 2^{j_i}-1}       E^{(j_{i+1}^{\rm cut})}_{\rm cut}         \right)  \notag \\
		&+\sum_{i = 1}^m    \left(   \sum_{j_1=0}^{j^{\rm cut}_{1}} \sum_{k_1=0}^{N_1 2^{j_1}-1}   \dots     \sum_{j_i=0}^{j^{\rm cut}_{i}} \sum_{k_i=0}^{N_i 2^{j_i}-1}       E_{\rm well}^{(j_i,k_i)}          \right)   \notag \\
		&+\sum_{i = 1}^m    \left(   \sum_{j_1=0}^{j^{\rm cut}_{1}} \sum_{k_1=0}^{N_1 2^{j_1}-1}    \dots    \sum_{j_i=0}^{j^{\rm cut}_{i}} \sum_{k_i=0}^{N_i 2^{j_i}-1}   \varepsilon    E_{\rm surf}^{(j_i,k_i)}          \right)  \notag \\
		&+    \sum_{j_1=0}^{j^{\rm cut}_{1}} \sum_{k_1=0}^{N_1 2^{j_1}-1}    \dots    \sum_{j_m=0}^{j^{\rm cut}_{m}} \sum_{k_m=0}^{N_m 2^{j_m}-1} E_{\rm pay}^{(j_m,k_m)} 
	\end{align}
with the convention $j^{\rm cut}_{0} = 0$ and $N_0 2^{j_0}-1= 0$.
Considering the $i$th branching construction, in the spirit of Remark~\ref{rem:firstorderenergycont}, $E^{(j_{i}^{\rm cut})}_{\rm cut}$ collects the contribution of the cut-off part, and $E_{\rm well}^{(j_i,k_i)}$ corresponds to the elastic energy where the gradient is close to the wells,
$E_{\rm surf}^{(j_i,k_i)}$ denotes the surface energy  
and $ E_{\rm pay}^{(j_i,k_i)} $ denotes the part, where the gradient is far away from $K$ and closer to an auxiliary matrix $J_j$, $j = 1,\dots ,4$.
\end{lemma}
\begin{proof}
Recalling $E^{(1\rm st)}$ in Remark~\ref{rem:firstorderenergycont}, which is the total energy contribution of the first-order branching construction, it holds that
\begin{align*}
	E^{(1\rm st)} &= E^{(j_1^{\rm cut})} + \sum_{j_1=0}^{j^{\rm cut}_{1}} \sum_{k_1=0}^{N_1 2^{j_1}-1}  \left(E_{\rm well}^{(j_1,k_1)} + E_{\rm pay}^{(j_1,k_1)} + \varepsilon E_{\rm surf}^{(j_1,k_1)}  \right) .
\end{align*}
	For the second-order construction, the energy $E_{\rm pay}^{(j_1,k_1)}$ changes and is rewritten via $$ E_{\rm pay}^{(j_1,k_1)} =   E^{(j_2^{\rm cut}  )} + \sum_{j_2=0}^{j^{\rm cut}_{2}} \sum_{k_2=0}^{N_2 2^{j_2}}   \left( E_{\rm well}^{(j_2 ,k_2 )} + E_{\rm pay}^{(j_2 ,k_2 )} + \varepsilon E_{\rm surf}^{(j_2 ,k_2 )} \right) .$$
More precisely, given $E_{\rm pay}^{(j_{m-1},k_{m-1})}$, we find
	$$ E_{\rm pay}^{(j_{m-1},k_{m-1})} =   E^{(j_m^{\rm cut}  )} + \sum_{j_m=0}^{j^{\rm cut}_{m}} \sum_{k_m=0}^{N_m 2^{j_m}-1}   \left( E_{\rm well}^{(j_m,k_m)} + E_{\rm pay}^{(j_m,k_m)} + \varepsilon E_{\rm surf}^{(j_m,k_m)} \right) .$$
Thus, an induction shows \eqref{energycontribution}.
\end{proof}
In the next lemma, we deduce bounds for the terms in \eqref{energycontribution}.
\begin{lemma}[$m$th-order construction]\label{lem:iteratedformula}
Consider the construction of Lemma~\ref{lem:preliminaryiteration}.
Then, we have
\begin{align}\label{Controlpaypart}
  \sum_{j_1=0}^{j^{\rm cut}_{1}} \sum_{k_1=0}^{N_1 2^{j_1}-1}    \dots    \sum_{j_m=0}^{j^{\rm cut}_{m}} \sum_{k_m=0}^{N_m 2^{j_m}-1} E_{\rm pay}^{(j_m,k_m)}  \lesssim   \frac{1}{2^m},
\end{align}
as well as
\begin{align}\label{boundEwell}
	 & \quad \ \sum_{i = 1}^m    \left(   \sum_{j_1=0}^{j^{\rm cut}_{1}} \sum_{k_1=0}^{N_1 2^{j_1}-1}   \dots     \sum_{j_i=0}^{j^{\rm cut}_{i}} \sum_{k_i=0}^{N_i 2^{j_i}-1}       E_{\rm well}^{(j_i,k_i)}          \right) \notag \\
	 &\lesssim  \sum_{i = 1}^m i^2 (1-\theta)^i 2^{-i} \sum_{j_1=0}^{j^{\rm cut}_{1}}   \theta^{j_1} \dots \sum_{j_{i-1}=0}^{j^{\rm cut}_{i-1}}   \theta^{j_{i-1}}  \sum_{j_i=0}^{j^{\rm cut}_{i}}   (\theta^{(1-p)} 2^{-p} )^{j_i}      r_i^p r_{i-1}^{-p}  ,
\end{align}
and  
\begin{align}\label{boundEsurf}
& \quad \ \sum_{i = 1}^m  \left(   \sum_{j_1=0}^{j^{\rm cut}_{1}} \sum_{k_1=0}^{N_1 2^{j_1}-1}    \dots    \sum_{j_i=0}^{j^{\rm cut}_{i}} \sum_{k_i=0}^{N_i 2^{j_i}-1}   \varepsilon    E_{\rm surf}^{(j_i,k_i)}          \right)  \notag \\
	& \lesssim \varepsilon 	\sum_{i = 1}^m  i^2  (2-2\theta)^{i/2}  \left(   \sum_{j_1=0}^{j^{\rm cut}_{1}}    (2\theta)^{j_{1}}  \dots    \sum_{j_i=0}^{j^{\rm cut}_{i}} (2\theta)^{j_{i}}   \right)           r_{i}^{-1} .
\end{align}
Moreover,
\begin{align}\label{boundEcut}
& \quad \ \sum_{i = 0}^{m-1}  \left(   \sum_{j_0=0}^{j^{\rm cut}_{0}} \sum_{k_0=0}^{N_0 2^{j_0}-1}   \dots     \sum_{j_i=0}^{j^{\rm cut}_{i}} \sum_{k_i=0}^{N_i 2^{j_i}-1}       E^{(j_{i+1}^{\rm cut})}_{\rm cut}         \right) \notag \\
 & \lesssim   \sum_{i = 1}^m i^2 (1-\theta)^i 2^{-i}  \theta^{j^{\rm cut}_{1}+1} \cdots  \theta^{j^{\rm cut}_{i-1}+1}    (\theta^{(1-p)} 2^{-p} )^{j^{\rm cut}_{i}+1}      r_i^p r_{i-1}^{-p}     \notag \\
& \qquad + \varepsilon 	\sum_{i = 1}^m  i^2  (2-2\theta)^{i/2}  \left(      (2\theta)^{j^{\rm cut}_{1}+1}  \cdots      (2\theta)^{j^{\rm cut}_{i}+1}   \right)           r_{i}^{-1}.
\end{align}
\end{lemma}
\begin{proof}
The parts, where we pay lots of elastic energy have been reduced by the factor $\frac{1}{2}$ after the first iteration, and we find $\sum_{j_2=0}^{j^{\rm cut}_{2}} \sum_{k_2=0}^{N_2 2^{j_2}} E_{\rm pay}^{(j_2,k_2)} \leq 
\frac{1}{2} E_{\rm pay}^{(j_1,k_1)}$ after the first replacement. As shown in Remark~\ref{rem:firstorderenergycont}, $\sum_{j_1=0}^{j^{\rm cut}_{1}} \sum_{k_1=0}^{N_1 2^{j_1}}  E_{\rm pay}^{(j_1,k_1)}$ is uniformly bounded, which shows
\eqref{Controlpaypart} by an iteration.
In order to show \eqref{boundEwell} and \eqref{boundEsurf}, we derive a bound for every addend $i \in \{1, \dots ,m\}$ in these inequalities.
As a preliminary step, we prove formulas for $\ell_{j_m}$ and $h_{j_m}$.
In view of \eqref{eq:lengthdef} and \eqref{def:lengthsblocksit}, and induction shows
\begin{align*}
	\ell_{j_m} = \frac{(1-\theta)^{k-1}}{8^{k-1}} \frac{\prod_{i = 1}^{k-1}\theta^{j_{2i}}}{\prod_{i=1}^{k} 2^{j_{2i-1}} N_{2i-1} }  \text{ and } h_{j_m} = \frac{1}{2}  \frac{(1-\theta)^{k}}{8^{k-1}} \frac{\prod_{i=1}^{k} \theta^{j_{2i-1}}}{\prod_{i = 1}^{k-1} 2^{j_{2i}} N_{2i} }
\end{align*}
for $m = 2k -1$, $k \in \N$, and
\begin{align*}
	\ell_{j_m} = \frac{1}{2} \frac{(1-\theta)^{k}}{8^{k-1}} \frac{\prod_{i = 1}^k\theta^{j_{2i-1}}}{\prod_{i=1}^{k} 2^{j_{2i}} N_{2i} }  \text{ and } h_{j_m} =   \frac{(1-\theta)^{k}}{8^{k}} \frac{\prod_{i=1}^{k} \theta^{j_{2i}}}{\prod_{i = 1}^{k} 2^{j_{2i-1}} N_{2i-1} }
\end{align*}
for $m = 2k $, $k \in \N$.
 We now set $N_m = \lceil (2\theta)^{j_{m-1}} r_m^{-1} r_{m-2} \rceil$, where $j_0 = 0$, $r_{-1} = r_0 = 1$. 
Let us now prove that
\begin{align}\label{bound:Nm}
	(2\theta)^{j_{m-1}} r_m^{-1} r_{m-2} \lesssim N_m \lesssim m (2\theta)^{j_{m-1}} r_m^{-1} r_{m-2}
\end{align}
for $\varepsilon$ sufficiently small.
Let us recall that $j_m^{\rm cut}$ is the largest index, for which $\ell_{j_m} \leq h_{j_m}$.
By definition of $\ell_{j_m^{\rm cut}+1}$ and $h_{j_m^{\rm cut} +1}$, the fraction
$\ell_{j_{m}} h_{j_{m}}^{-1} $ is uniformly bounded for all $j_{m} = 0, \dots,j_m^{\rm cut} +1$. 
Thus, the formulas for $\ell_{j_m}$ and $h_{j_m}$ imply that there is a positive constant $C>0$ such that
\begin{align*}
	\prod_{i=1}^m \left( (2\theta)^{j_i} N_i \right)^{(-1)^{m-i}} \geq C.
\end{align*}
This yields
\begin{align}\label{eq:equationnowuse}
	\frac{r_{m+1}}{r_{m}} \prod_{i=1}^m \left(\frac{ N_i r_i }{  (2\theta)^{j_{i-1}} r_{i-2}  } \right)^{(-1)^{m-i}}   = \frac{r_{m+1}}{(2\theta)^{j_m} r_{m-1}}\prod_{i=1}^m \left( (2\theta)^{j_i} N_i \right)^{(-1)^{m-i}} \geq C \frac{r_{m+1}}{(2\theta)^{j_m} r_{m-1}}.
\end{align}
Further, it holds that
\begin{align}\label{eq:forconclusion}
	1\leq \frac{ N_i r_i }{  (2\theta)^{j_{i-1}} r_{i-2}  } \leq \frac{ r_i }{  (2\theta)^{j_{i-1}} r_{i-2}  } + 1
\end{align}
for all $i= 1, \dots, m$.
The previous upper bound for $i = m$, \eqref{eq:equationnowuse} for $m-1$ instead of $m$, the previous lower bound for $i = m-1$, and $r_{m} r_{m-1}^{-1} \leq 1$ yield
\begin{align*}
	& \quad \ \prod_{i=1}^{m} \left(\frac{ N_i r_i }{  (2\theta)^{j_{i-1}} r_{i-2}  } \right)^{(-1)^{m-i}} \\
	&\leq \prod_{i=1}^{m-1} \left(\frac{ N_i r_i }{  (2\theta)^{j_{i-1}} r_{i-2}  } \right)^{(-1)^{m-i}} + C \frac{r_{m}}{(2\theta)^{j_{m-1}} r_{m-2}}  \prod_{i=1}^{m-1} \left(\frac{ N_i r_i }{  (2\theta)^{j_{i-1}} r_{i-2}  } \right)^{(-1)^{m-i}}   \\
	&\leq \left(\frac{ N_{m-1} r_{m-1} }{  (2\theta)^{j_{m-2}} r_{m-3}  } \right)^{-1}  \prod_{i=1}^{m-2} \left(\frac{ N_i r_i }{  (2\theta)^{j_{i-1}} r_{i-2}  } \right)^{(-1)^{m-i}} + C \frac{r_{m}}{r_{m-1}} \\
		&\leq  \prod_{i=1}^{m-2} \left(\frac{ N_i r_i }{  (2\theta)^{j_{i-1}} r_{i-2}  } \right)^{(-1)^{m-i}}  + C .
\end{align*}
An iteration of the latter inequality gives
\begin{align*}
	 \prod_{i=1}^{m} \left(\frac{ N_i r_i }{  (2\theta)^{j_{i-1}} r_{i-2}  } \right)^{(-1)^{m-i}} \leq C m.
\end{align*}
In view of \eqref{eq:equationnowuse}, we hence find that
\begin{align*}
	\frac{r_{m+1}}{(2\theta)^{j_m} r_{m-1}} \leq C m \frac{r_{m+1}}{r_{m}}.
\end{align*}
\eqref{eq:forconclusion} and $\frac{r_{m+1}}{r_{m}}\leq 1$ eventually show \eqref{bound:Nm}.
Using \eqref{bound:Nm} and the formulas for $\ell_{j_m}$ and $h_{j_m}$, we obtain the lower bounds
\begin{align}\label{iterationformulaellh}
	\ell_{j_m} \gtrsim m^2 \left( \frac{1-\theta}{8} \right)^{m/2} r_m \prod_{s=1}^m 2^{-j_s}  \quad \text{ and } \quad h_{j_m} \gtrsim  m^2 \left( \frac{1-\theta}{8} \right)^{m/2} \theta^{j_m} r_{m-1} \prod_{s=1}^{m-1} 2^{-j_s}
\end{align}
and upper bounds
\begin{align}\label{iterationformulaellh2}
	\ell_{j_m} \lesssim \left( \frac{1-\theta}{8} \right)^{m/2} r_m \prod_{s=1}^m 2^{-j_s}  \quad \text{ and } \quad h_{j_m} \lesssim  \left( \frac{1-\theta}{8} \right)^{m/2} \theta^{j_m} r_{m-1} \prod_{s=1}^{m-1} 2^{-j_s}.
\end{align}
 We now prove \eqref{boundEwell}--\eqref{boundEcut}.
The self-similar construction gives
 \begin{align}\label{eq:elasticenergybound}
   \sum_{i = 1}^m    \left(   \sum_{j_1=0}^{j^{\rm cut}_{1}} \sum_{k_1=0}^{N_1 2^{j_1}-1}   \dots     \sum_{j_i=0}^{j^{\rm cut}_{i}} \sum_{k_i=0}^{N_i 2^{j_i}-1}       E_{\rm well}^{(j_i,k_i)}          \right)   \lesssim  \sum_{i = 1}^m    \left(   \sum_{j_1=0}^{j^{\rm cut}_{1}}  \dots \sum_{j_i=0}^{j^{\rm cut}_{i}}  \left(\prod_{s=1}^i N_s 2^{j_s} \right)        E_{\rm well}^{(j_i,k_i)}          \right) .
 \end{align}
As we have not indicated the branching patterns which are generated by reflection and shearing,
the energy contribution of the elastic energy
$E_{\rm well}^{(j_i,k_i)}$ on every block has to be controlled by \eqref{eq:elbuildingblock} multiplied with the factor $4^i$.
Let us provide a bound for the last summand $i = m$. Using \eqref{bound:Nm} once again, 
by \eqref{eq:elbuildingblock}, \eqref{iterationformulaellh}, and \eqref{iterationformulaellh2} it holds that
\begin{align*}
	 \left(\prod_{s=1}^m N_s 2^{j_s} \right)        E_{\rm well}^{(j_m,k_m)}  &  \lesssim 4^{m} \frac{\ell_{j_m}^{p+1}}{h_{j_m}^{p-1}} \prod_{s=1}^m N_s 2^{j_s} \\
		&\lesssim m^2 \left( \frac{1-\theta}{2} \right)^{m} r_m^{p} r_m^{-p} 2^{-p j_m}  \theta^{j_m(1-p)}     \left(\prod_{s=1}^{m-1} \theta^{j_s}\right)   .
\end{align*}
The previous estimate and \eqref{eq:elasticenergybound} conclude \eqref{boundEwell}.
In order to prove \eqref{boundEsurf}, we note that also
$E_{\rm surf}^{(j_i,k_i)}$ has to be controlled via the estimate \eqref{eq:elbuildingblock} multiplied with the factor $4^i$ as the surface energy which arises from reflection and shearing has to be counted.
We again only provide the bound for the last summand $i = m$. By \eqref{iterationformulaellh}, \eqref{iterationformulaellh2}, \eqref{bound:Nm}, and \eqref{eq:elbuildingblock} it holds that
\begin{align*}
  \left(\prod_{s=1}^m N_s 2^{j_s} \right)     E_{\rm surf}^{(j_m,k_m)} & \lesssim 4^m h_{j_m}   \left(\prod_{s=1}^m N_s 2^{j_s} \right) \\
  & \lesssim m^2 4^m \left(\left( \frac{1-\theta}{8} \right)^{m/2} \theta^{j_m} r_{m-1} \prod_{s=1}^{m-1} 2^{-j_s} \right)  \left(\prod_{s=1}^m (2\theta)^{j_{s-1}} r_s^{-1} r_{s-2} 2^{j_s} \right) \\
    & \lesssim m^2  (2-2\theta)^{m/2}   \left(\prod_{s=1}^m (2\theta)^{j_{s}}   \right) r_m^{-1} .
\end{align*}
This shows
\eqref{boundEsurf}. \eqref{boundEcut}  follows by \eqref{eq:cut-off}. Indeed, the computation follows along similar lines as the cut-off energy is controlled via \eqref{eq:cut-off} instead of \eqref{eq:elbuildingblock}. The bounds, however, are of the same order.
\end{proof}

\subsection{Proof of Theorem~\ref{mainth:upper}}\label{sec:updet}
We now determine the optimal length scales for the bounds derived in Lemma~\ref{lem:iteratedformula} and hence provide the proof of the upper bound.

\begin{proof}[Proof of Theorem~\ref{mainth:upper}]
Let us identify the optimal length scales for the construction in the previous section. We also refer to \cite[Lemma~3.2]{W97}, where the optimal length scale for an iterated laminate construction was determined.
More precisely, by Lemma~\ref{lem:preliminaryiteration} and Lemma~\ref{lem:iteratedformula}, there are competitors such that
\begin{align}\label{eq:optimizationproblemimproved}
	\Energyinf{p} \lesssim 2^{-m} +  r_1^p + \sum_{i=1}^m i^2 \left( \frac{1-\theta}{2} \right)^i \frac{r_{i+1}^p}{r_i^p}  + \varepsilon \sum_{i=1}^m i^2  (2-2\theta)^{i/2}  r_i^{-1},
\end{align}
provided that
$\theta \in(2^{-p/(p-1)},1/2)$ if $p>1$ and $\theta \in (0,1/2)$ if $p=1$, respectively.
Heuristically, the length scales $r_{i}$, $i \in \N$, should be chosen in such a way that every term in the sum above generates
the same asymptotic behaviour, leading to the preliminary choice
$r_1^p \sim \big( \frac{2}{1-\theta} \big)^i  \frac{r_{i+1}^p}{r_i^p}$ for all $i \in \N$.
By induction, one obtains the formula
	$$r_i \sim r_1^i \left( \left( \frac{2}{1-\theta} \right)^{\frac{i(i+1)}{4}} \right)^{2/p}.$$
This choice, however, generates the $m+1$-fold product of the term $r_1^p$ on the right-hand side of \eqref{eq:optimizationproblemimproved}
as we sum over $m$ terms which do not depend on $m$. Due to quantity $2^{-m}$ in \eqref{eq:optimizationproblemimproved}, we necessarily need $m \to \infty$ as $\varepsilon \to 0$ in order to obtain a bound converging to zero.
Recalling the geometric series, we hence make the ansatz
\begin{align}\label{eq:ansatzoptimalscaling}
	r_i = r_1^i \left( \alpha^{\frac{i(i+1)}{4}} \right)^{2/p} = r_1^i   \alpha^{\frac{i(i+1)}{2p}}  \text{ for } \alpha \in \Big(1, \frac{2}{1-\theta} \Big) \text{ fixed} .
\end{align}
Then, an elementary computation shows that
\begin{align*}
	i^2 \left( \frac{1-\theta}{2} \right)^i \frac{r_{i+1}^p}{r_i^p} &= i^2 \left( \frac{1-\theta}{2} \right)^i \frac{ \Big(r_1^{i+1} \alpha^{\frac{(i+1)(i+2)}{2p}} \Big)^p}{\Big(r_1^i \alpha^{\frac{i(i+1)}{2p}}\Big)^p} = \alpha i^2 \left( \frac{1-\theta}{2} \right)^i \alpha^i r_1^p.
\end{align*}
Thus, we obtain the estimate
\begin{align}\label{eq:optscale_elasticenergy}
	  r_1^p + \sum_{i=1}^m i^2  \left( \frac{1-\theta}{2} \right)^i \frac{r_{i+1}^p}{r_i^p} \leq   C r_1^p
\end{align}
for a positive constant $C >0$, depending only on $\alpha$ and $\theta$.  
Further, the terms $r_1^p$, $2^{-m}$, and $\varepsilon \sum_{i=1}^m i^2 (2-2\theta)^{i/2}   r_i^{-1}$ should capture the same scaling, see \eqref{eq:optimizationproblemimproved}. 
To this end, we make the ansatz
\begin{align}\label{choiceofr_1}
r_1 = \exp \Big(-  \eta \log(\varepsilon^{-1})^{1/2} \Big) \quad 
\text{ as well as } \quad 
m = \left\lceil \gamma  \log(\varepsilon^{-1})^{1/2}  \right\rceil
\end{align}
for some $\gamma$ and $\eta > 0$.
By this choice, we immediately see that
\begin{align}\label{eq:optscale_volume}
	r_1^p + 2^{-m} \leq   C \exp \Big(- \eta p \log(\varepsilon^{-1})^{1/2} \Big) + C \exp \Big(- \gamma \log(2) \log(\varepsilon^{-1})^{1/2} \Big).
\end{align} 
Next, we analyse the surface energy.
For $\beta > (2-2\theta)^{1/2}$, there is an index $i_0$ such $i^2 (2-2\theta)^{i/2} < \beta^i$ for $i \geq i_0$. Thus, an index change and 
\eqref{eq:ansatzoptimalscaling} show that the surface energy can be rewritten as
\begin{align}\label{eq:surfacenergyest}
\varepsilon \sum_{i=1}^{m} i^2 (2-2\theta)^{i/2} r_i^{-1}   & = 2^{-m} \varepsilon \beta^m 2^m \sum_{i=0}^{m-1} \beta^{-i} r_{m-i}^{-1}  = 2^{-m} \varepsilon m^2 \beta^m 2^m \sum_{i=0}^{m-1}  \beta^{-i}  r_1^{i-m}   \alpha^{\frac{-(m-i)(m-i+1)}{2p}}   \notag  \\
& = 2^{-m} \varepsilon \beta^m 2^m \alpha^{\frac{-m^2-m}{2p}}  r_1^{-m} \sum_{i=0}^{m-1}  \alpha^{\frac{i}{2p}}  \beta^{-i}  \alpha^{\frac{-i^2}{2p}}    r_1^{i} \alpha^{\frac{mi}{p}} \notag   \\
& = 2^{-m} A(m) \sum_{i=0}^{m-1}  B(i)    D(i,m)  ,
\end{align}
where
\begin{align*}
A(m) \defas \varepsilon \beta^m 2^m \alpha^{\frac{-m^2-m}{2p}}  r_1^{-m}, \quad B(i) \defas \alpha^{\frac{i}{2p}}  \beta^{-i}  \alpha^{\frac{-i^2}{2p}}, \quad D(i,m) \defas r_1^{i} \alpha^{\frac{mi}{p}}.
\end{align*}
As $2^{-m}$ already captures the desired scaling, see \eqref{eq:optscale_volume},
it suffices to prove that
$A(m)$ and $D(i,m)$ are uniformly bounded and that $B(i) < q^i$ for some $q \in (0,1)$.
The latter condition is satisfied for $\alpha > 1$ as $i$ appears quadratically in the exponent.
Using \eqref{choiceofr_1}, we have
\begin{align*}
	A(m) \leq C \exp \bigg( \Big(\gamma\eta -1 - \frac{\log(\alpha)}{2p} \Big) \log(\varepsilon^{-1}) \bigg) \exp \bigg( \Big(\gamma \log(2) + \gamma \log(\beta) - \gamma \frac{\log(\alpha) } {2p}   \Big) \log(\varepsilon^{-1})^{1/2} \bigg).
\end{align*}
Hence, boundedness is satisfied for $\varepsilon$ small, provided that $$\gamma \eta < 1 +  \frac{\gamma^2\log(\alpha)}{2p}.$$
Further, if $\log(\alpha) \gamma p^{-1} \leq \eta$, it holds that
\begin{align*}
	D(i,m) \leq C \exp \big( ( \log(\alpha) \gamma p^{-1} - \eta ) \log(\varepsilon^{-1})^{1/2} i  \big) \leq C.
\end{align*}
In order to obtain a sharp bound, we have to choose $C_{\rm upp}$ as large as possible.
In view of \eqref{eq:optimizationproblemimproved}, \eqref{eq:optscale_elasticenergy}, and \eqref{eq:optscale_volume} we thus have
\begin{align*}
	2^{-m} +  r_1^p + \sum_{i=1}^m i^2 \left( \frac{1-\theta}{2} \right)^i \frac{r_{i+1}^p}{r_i^p}  + \varepsilon \sum_{i=1}^m i^2 (2-2\theta)^{i/2}  r_i^{-1} \leq \exp \big( - C_{\rm upp}  \log(\varepsilon^{-1})^{1/2}\big),
\end{align*}
where
\begin{align*}
	C_{\rm upp} = \min \{ \eta p , \gamma \log(2)      \} \quad \text{ with } \quad \log(\alpha) \gamma p^{-1} \leq \eta, \ \gamma \eta < 1 +  \frac{\gamma^2 \log(\alpha)}{2p}, \text{ and }  \alpha \in \Big(1, \frac{2}{1-\theta} \Big).
\end{align*}
An optimization, see Remark~\ref{rem:optimalityupper} below, leads to the choice $\alpha = 2$, $\eta = \sqrt{2 p^{-1}} \log(2)^{1/2}$, $\gamma = t \sqrt{2p} \log(2)^{-1/2}$ for $t \in (0,1)$.
The parameter $t \in (0,1)$ is needed in order to ensure that $\eta \gamma - \gamma^2 \log(2) (2p)^{-1} = 2t -t^2 < 1$. Further, $\eta > \gamma \log(2) p^{-1}$.
By sending $t \to 1$, we have 
\begin{align*}
	C_{\rm upp} \to \sqrt{2p \log(2)}.
\end{align*}
This concludes the proof.
\end{proof}

\begin{remark}[Comparison with the optimization in \cite{W97}]\label{rem:optimalityupper}
Let us compare our bounds with the analogous ones from \cite{W97} for the laminate construction.
Therein, the bounds are similar to \eqref{eq:optimizationproblemimproved}, but $r_{i+1}^p r_i^{-p}$ and $r_1^p$ need to be replaced by $r_{i+1} r_i^{-1}$ and $r_1$, respectively.
Further, both prefactors $i^2 ( \frac{1-\theta}{2} )^i $ and $ i^2 (2-2\theta)^{i/2} $ in the elastic and surface contribution are substituted by $2^{-i}$.
The exponent $p$ arises from branching construction and is well-known in finite-order constructions.
The other two prefactors also occur in branching constructions, see \cite{KM2} for the case of first-order data,
but were not tracked for higher-order data in \cite{RT23} as these do not influence the scaling in $\varepsilon$ if the iteration parameter $m$ is fixed.
In our context, as we also optimize in $m$, the parameters contribute in \eqref{eq:optimizationproblemimproved}. Hence, one rigorously needs to keep track of these factors because they might influence the scaling.

Let us first comment on $i^2 (2-2\theta)^{i/2}$, which is related to $\beta^i$ in the previous proof. In fact, the computation below \eqref{eq:surfacenergyest} imply that the terms remain bounded for any $\beta > 0$.
For the discussion of the prefactor $i^2  ( \frac{1-\theta}{2}  )^i $, let us discuss that the choice $\alpha = 2$ is optimal in the previous proof. The constraint
 \begin{align*}
	\eta \in \left[ \log(\alpha) \gamma p^{-1}     ,    \gamma^{-1} +  \frac{\gamma \log(\alpha)}{2p} \right)
\end{align*}
needs that $\gamma^2 < \frac{2p}{\log(\alpha)}$, otherwise, the interval is empty. For a fixed $\gamma$, $C_{\rm upp}$ becomes large if $\eta$ approaches the upper bound of the interval, i.e.,
$$C_{\rm upp} \to \min \left\{ \gamma^{-1} p +  \frac{\gamma \log(\alpha)}{2} , \gamma \log(2)      \right\}.$$
Balancing the minimum problem leads to the choice
 $ \gamma = p^{1/2} \big( \log(2)-\frac{ \log(\alpha)}{2} \big)^{-1/2} $, implying that
 $$C_{\rm upp} \to \log(2) p^{1/2} \bigg( \log(2)-\frac{ \log(\alpha)}{2} \bigg)^{-1/2}.$$
For $\alpha >2$, the previous choice of $\gamma$ is not admissible due to $\gamma^2 = \frac{p}{\log(2) - \log(\alpha) /2} > \frac{2p}{\log(\alpha)}$. This contradicts the constraint that the interval has a positive length.
In this case, the choice $\gamma \to (2p)^{1/2}\log(\alpha)^{-1/2}$ leads to
 $$C_{\rm upp} \to \log(2) (2p)^{1/2} \log(\alpha)^{-1/2},$$
showing that $\alpha = 2$ is optimal. Hence, the additional factor $(1-\theta)^i$ in \eqref{eq:optimizationproblemimproved} is not exploited to obtain a better scaling behaviour.
As a summary, the exponent $p$ improves the scaling behaviour, but the additional factors in $\theta$ do not.
\end{remark}

\subsection*{Acknowledgements}
The author gratefully thanks Antonio Tribuzio and Angkana Rüland for helpful discussions.
The work is funded by the Deutsche Forschungsgemeinschaft (DFG,
German Research Foundation) -- CRC 1720 -- 539309657 and by the Deutsche Forschungsgemeinschaft (DFG, German Research Foundation)
under Germany's Excellence Strategy -- EXC-2047/1 -- 390685813.

\subsection*{AI Statement}
The mathematical arguments and a complete draft of this paper were developed by the author without generative AI. 
Subsequent use of ChatGPT was limited to the assistance of algebraic manipulations in \eqref{eq:assistchat}, \eqref{bound:Nm}, computations in Remark~\ref{rem:optimalityupper},
and the creation of Figure~\ref{fig:branching}.
All the proofs contained in this manuscript were written by the author, who takes full responsibility for the correctness of the statements. This article does not contain AI-written text.

\subsection*{Data Availability Statement}
No data has been used in this work.

\subsection*{Conflict of interest}
 The author has no conflict of interest to declare.

\printbibliography

\end{document}